\documentclass[final,3p,times]{elsarticle}

\usepackage{graphicx,epstopdf,amsmath,amssymb,txfonts}
\usepackage{algorithm}
\usepackage{algorithmic}
\usepackage{subfigure,color}
\usepackage[colorlinks]{hyperref}
\usepackage{hypernat}
\usepackage{multirow}
\numberwithin{equation}{section}
\graphicspath{{Figs/}}
\usepackage{booktabs}
\usepackage{listings}
\usepackage{lineno,hyperref}

\biboptions{sort&compress}

\makeatletter

\newcommand{\Rmnum}[1]{\expandafter\@slowromancap\romannumeral #1@}
\makeatother

\journal{Computers \& Mathematics with Applications}

\newtheorem{theorem}{\textbf{Theorem}}
\newtheorem{Algorithm}{\textbf{Algorithm}}
\newtheorem{remark}{\textbf{Remark}}
\newtheorem{example}{\textbf{Example}}
\newenvironment{Proof}[1][Proof.]{\begin{trivlist}
\item[\hskip \labelsep {\bfseries #1}]}{\ep\end{trivlist}}
\newcommand{\ep}{\hfill\rule{0.15cm}{0.35cm}\vskip 0.3cm}

\begin{document}

\begin{frontmatter}



\title{A {\color{red}multigrid-reduction-in-time} solver with a new two-level convergence for unsteady fractional Laplacian problems}


\author[myadda]{Xiaoqiang Yue}
\ead{yuexq@xtu.edu.cn}

\author[myaddc]{Kejia Pan\corref{mycoaur}}
\ead{kejiapan@csu.edu.cn}

\author[myadda]{Jie Zhou}
\ead{zhouj@xtu.edu.cn}

\author[myaddb]{Zhifeng Weng}
\ead{zfwmath@163.com}

\author[myadda]{Shi Shu\corref{mycoraur}}
\ead{shushi@xtu.edu.cn}

\author[myadda]{Juan Tang}

\address[myadda]{Key Laboratory of Intelligent Computing \& Information Processing of Ministry of Education,
Hunan Key Laboratory for Computation and Simulation in Science and Engineering,
School of Mathematics and Computational Science,
Xiangtan University, Xiangtan 411105, China}
\address[myaddc]{School of Mathematics and Statistics, HNP-LAMA,
Central South University, Changsha 410083, China}
\address[myaddb]{School of Mathematical Sciences, Fujian Province University Key Laboratory of
Computation Science,\\ Huaqiao University, Quanzhou 362021, China}

\cortext[mycoaur]{Xiaoqiang Yue and Kejia Pan contributed equally to this work and should be considered co-first authors.}
\cortext[mycoraur]{Corresponding author.}

\begin{abstract}
The multigrid-reduction-in-time (MGRIT) technique has proven to be successful
in achieving higher run-time speedup by exploiting parallelism in time.
The goal of this article is to develop and analyze a MGRIT algorithm,
using FCF-relaxation with time-dependent time-grid propagators,
to seek the finite element approximations of unsteady fractional Laplacian problems.
{\color{red}The multigrid with line smoother proposed in [L. Chen, R. H. Nochetto, et al.,
Math. Comp. 85 (2016) 2583--2607] is chosen to be the spatial solver.}
Motivated by [B. S. Southworth, SIAM J. Matrix Anal. Appl. 40 (2019) 564--608],
we provide a new temporal eigenvalue approximation property and then
deduce a generalized two-level convergence theory which removes the previous unitary diagonalization assumption
on the fine and coarse time-grid propagators
required in [X. Q. Yue, S. Shu, et al.,
Comput. Math. Appl. 78 (2019) 3471--3484].
Numerical computations are included to confirm the theoretical predictions and demonstrate
the sharpness of the derived convergence upper bound.
\end{abstract}

\begin{keyword}
fractional Laplacian, Caffarelli-Silvestre extension,
reduction-based multigrid, multigrid-in-time,
finite element, anisotropic mesh.
\end{keyword}

\end{frontmatter}

\section{Introduction}

Time-dependent problems involving fractional Laplacian operator are of particular interest in fractional calculus \cite{p-01,w-02,m-02}.
The fractional Laplacian, which is the infinitesimal generator of a standard isotropic $\alpha$-stable L\'{e}vy flight \cite{s-01},
has been used in a broad range of applications
including quantum mechanics \cite{l-03}, finance \cite{c-06},
elasticity \cite{d-01}, porous medium \cite{p-02} and stochastic dynamics \cite{g-04}.
Further details on this topic and equivalent characterizations can be found in the review article \cite{l-01}.
We are concerned in this paper with a non-intrusive parallel time integration with multigrid to
a class of unsteady fractional Laplacian problems of the form
\begin{align}
  \label{eq-01}
  \left \{
    \begin{aligned}
    & \partial_tu({\bf x},t)=-(-\Delta)^{\alpha/2}u({\bf x},t)+f({\bf x},t), & t\in I:=(0,T],~{\bf x}\in \Omega, \\
    & u({\bf x},t)=0, & t\in I,~{\bf x}\in \partial\Omega, \\
    & u({\bf x},0)=\psi_0({\bf x}), & {\bf x}\in \Omega,
    \end{aligned}
  \right.
\end{align}
where the fractional order $\alpha\in(0,2)$,
solution domain {\color{red}$\Omega\subset \mathbb{R}^d$ ($d\ge1$)}
is an open and bounded subset with Lipschitz continuous boundary $\partial\Omega$.

One of the main challenges of problem \eqref{eq-01} lies in the nonlocality of the fractional Laplacian operator $(-\Delta)^{\alpha/2}$.
Numerous numerical approximations have been proposed to treat
the fractional Laplacian on bounded domains in recent years,
such as finite difference \cite{h-01,s-02,w-08,d-03,h-04,f-03,w-05,w-07,l-06,w-09},
finite element (FE) \cite{k-01,w-10,n-01,g-05,c-01,a-02,a-05,j-01},
finite volume \cite{c-07,g-08}, discontinuous Galerkin \cite{x-02,l-05,z-05,a-04}
and spectral (element) \cite{b-01,m-04,d-04,z-01,h-02,w-06,z-02,w-11,t-02,s-04} methods.
Utilizing a similar technique to the fractional Laplacian in \cite{c-05},
we can obtain the Caffarelli-Silvestre extension,
which is reformulated the original problem posed over $\Omega$
on the semi-infinite cylinder $\mathcal{C}=\Omega\times(0,+\infty)$ with dynamic boundary condition:
\begin{align}
  \label{eq-02}
  \left \{
    \begin{aligned}
    & \textrm{div}(z^\beta \nabla \mathcal{U}) = 0, & t\in I,~({\bf x},z)\in \mathcal{C}, \\
    & \mathcal{U}({\bf x},0,0) = \psi_0({\bf x}), & {\bf x}\in \Omega, \\
    & \mathcal{U}({\bf x},z,t)=0, & t\in I,~({\bf x},z)\in \partial_L\mathcal{C}, \\
    & d_\alpha\partial_t\mathcal{U}+\partial_\nu^\beta\mathcal{U}=d_\alpha f, & t\in I,~{\bf x}\in \Omega,~z=0,
    \end{aligned}
  \right.
\end{align}
where the parameter $\beta=1-\alpha$, the lateral boundary $\partial_L\mathcal{C}=\partial\Omega\times[0,+\infty)$,
the positive normalization constant $d_\alpha$ and the conormal exterior $\partial_\nu^\beta\mathcal{U}$ are respectively defined by
\begin{eqnarray*}
  d_\alpha = \frac{2^\beta\Gamma(1-\alpha/2)}{\Gamma(\alpha/2)}\quad\mbox{and}\quad
  \partial_\nu^\beta\mathcal{U} = -\lim_{z\rightarrow0^+}z^\beta \partial_z\mathcal{U}.
\end{eqnarray*}
A widely used approach to solve nonlocal problem \eqref{eq-01}
is based on the solution of local problem \eqref{eq-02},
however, incorporating one more dimension in space \cite{n-01,c-01}.
Specifically, by exploiting the intermediate function $\mathcal{U}({\color{red}{\bf x}},z,t): \mathcal{C}\times I\rightarrow \mathbb{R}$,
we immediately get the desired solution $u({\bf x},t): \Omega\times I\rightarrow \mathbb{R}$
subject to the expression $u({\bf x},t)=\mathcal{U}({\bf x},0,t)$ for all ${\bf x}\in\Omega$ and $t\in I$.
It is worthwhile to point out that the anisotropic spatial mesh
must be used to capture the singular behavior of $\mathcal{U}({\bf x},z,t)$ as $z\rightarrow0$ \cite{d-05},
and Nochetto et al. have proved the existence and uniqueness for solutions
of problems \eqref{eq-01} and \eqref{eq-02} (see \cite[Theorem 6]{n-02}) in a purely theoretical manner.

Moving from {\color{red}$d$}-dimensional problems to {\color{red}($d$+1)}-dimensional problems can be difficult for solvers.
Chen et al. developed uniformly convergent multilevel methods to tackle steady fractional Laplacian problems \cite{c-01},
which is an important constructive study in the literature.
However, a more severe challenge {\color{red}which we have to face} arises from the facts that (i)
the higher the spatial resolution signifies a comparative or greater scaling in the temporal direction,
(ii) future computing speed must rely on the increased concurrency provided by not faster but more processors.
The practical consequence is that solution algorithms limited to spatial parallelism
for problems with evolutionary behavior shall entail long overall computation time.
As a consequence, solution algorithms achieving parallelism in time have been of especially high demand.
Although time is sequential in essence, parallel-in-time methods are not new and 
may be dated back to \cite{n-03,m-05}.
A recent review of the literature can be found in \cite{g-07}.
At present, parareal in time \cite{l-02} and multigrid-reduction-in-time (MGRIT) \cite{f-02} are two of the most popular choices.
Wu and Zhou analyzed convergence properties of the parareal algorithm for space fractional diffusion problems \cite{w-01,w-04}
and time fractional equations with local time-integrators \cite{w-03} in constant time-steps.
Note that Gander and Vandewalle interpreted parareal as a two-level multigrid full approximation scheme \cite{g-02},
whereas Falgout et al. identified parareal with an optimal two-level MGRIT algorithm
and further pointed out that the obvious generalization of parareal to multiple levels couldn't produce an optimal algorithm \cite{f-02}.
In addition, the concurrency of parareal is still limited as its relatively large sequential coarse-grid solve.
Conversely, MGRIT, a truly multilevel algorithm implemented in the open-source package XBraid \cite{x-01},
has been proven to be rather effective and analyzed in two-level \cite{d-02} and multilevel \cite{h-03} settings,
for integer order basic parabolic and hyperbolic problems using L-stable Runge-Kutta schemes,
with the limitation on analysis that they only consider the same fine time-grid propagators.
Extensions on the two-level convergence analysis to the case of nonuniform temporal meshes are provided in \cite{y-01,s-03}.
However, we should point out that only MGRIT using F-relaxation was abstractly discussed
in \cite{s-03} and the unitary diagonalization assumption must be enforced in \cite{y-01}.
Moreover, from the survey of references, there are no calculations taking into account of MGRIT
to the fully implicit FE discretization of problem \eqref{eq-01}.

The arrangement of this work proceeds as follows.
Section \ref{sec2} introduces the functional framework
and outlines the fully implicit discrete FE scheme.
In Section \ref{sec3}, we present the MGRIT V(1,0)-cycle algorithm using FCF-relaxation
followed by a generalized two-level convergence analysis without the assumption of unitary diagonalization.
Section \ref{sec4} conducts numerical experiments to
illustrate the optimal order of convergence and the computational error-norms decay rate of the fully discrete FE scheme,
the robustness of Chen's multigrid with line smoother with respect to the fractional order $\alpha$
and the spatial mesh size, as well as the sharpness of the derived convergence upper bound on the two-level MGRIT algorithm.
Finally, in Section \ref{sec5}, we draw some conclusions and discuss future work.

\section{Mathematical preliminaries}\label{sec2}

This section is devoted to introducing some mathematical preliminaries,
including the variational formulation of the local problem \eqref{eq-02}
and its fully implicit discrete FE scheme.

\subsection{The variational formulation}

Let {\color{red}$\mathcal{D}\subset\mathbb{R}^{d+1}$} be open, $L^2(z^\beta,\mathcal{D})$ and $H^1(z^\beta,\mathcal{D})$
denote the weighted Lebesgue and Sobolev spaces \cite{g-06}
\begin{eqnarray*}
  L^2(z^\beta,\mathcal{D}) = \{w: \int_\mathcal{D} z^\beta w^2 < +\infty\}, 
  \quad H^1(z^\beta,\mathcal{D}) = \{v\in L^2(z^\beta,\mathcal{D}): |\nabla v|\in L^2(z^\beta,\mathcal{D})\},
\end{eqnarray*}
with norms
\begin{eqnarray*}
  \|w\|_{L^2(z^\beta,\mathcal{D})} = \Big{(}\int_\mathcal{D} z^\beta w^2\Big{)}^{\frac{1}{2}},
  \quad \|v\|_{H^1(z^\beta,\mathcal{D})} = \Big{(}\|v\|^2_{L^2(z^\beta,\mathcal{D})}
  + \|\nabla v\|^2_{L^2(z^\beta,\mathcal{D})}\Big{)}^{\frac{1}{2}},
\end{eqnarray*}
where $\nabla v$ is the distributional gradient of the measurable function $v$.

Since $\beta\in(-1,1)$, weight $z^\beta$ belongs to the Muckenhoupt class {\color{red}$A_2(\mathbb{R}^{d+1})$} \cite{m-03},
which implies that $H^1(z^\beta,\mathcal{D})$ is Hilbert.
In order to settle problem \eqref{eq-02}, we define
\begin{eqnarray*}
  \mathring{H}^1_L(z^\beta,\mathcal{D}) := \{v\in H^1(z^\beta,\mathcal{D}): v = 0 ~\mbox{on} ~\partial_L \mathcal{D}\},
  ~\mathbb{H}^{\alpha/2}(\Omega) = [H_0^1(\Omega),L^2(\Omega)]_{1-\alpha/2}~\mbox{and}~
  \mathbb{H}^{-\alpha/2}(\Omega)~\mbox{as the dual of}~\mathbb{H}^{\alpha/2}(\Omega).
\end{eqnarray*}
Thus, given a forcing term $f\in L^2(0,T;\mathbb{H}^{-\alpha/2}(\Omega))$ and an initial datum $\psi_0\in L^2(\Omega)$,
the function $u\in \mathbb{W}$ solves problem \eqref{eq-01}
if and only if its harmonic extension $\mathcal{U}\in \mathbb{V}$ solves problem \eqref{eq-02}, where
\begin{eqnarray*}
  \mathbb{W} := \{w\in L^\infty(0,T;L^2(\Omega))\cap L^2(0,T;\mathbb{H}^{\alpha/2}(\Omega)):\partial_tw\in L^2(0,T;\mathbb{H}^{-\alpha/2}(\Omega))\}
\end{eqnarray*}
and
\begin{eqnarray*}
  \mathbb{V} := \{v\in L^2(0,T;\mathring{H}^1_L(z^\beta,\mathcal{C})):\partial_t\mbox{tr}_\Omega v\in L^2(0,T;\mathbb{H}^{-\alpha/2}(\Omega))\}.
\end{eqnarray*}
Recall that the trace operator $\mbox{tr}_\Omega$ onto $\Omega\times\{z=0\}$ satisfies (see \cite[Proposition 2.5]{n-01})
\begin{eqnarray*}
  \mathring{H}^1_L(z^\beta,\mathcal{C})\ni w \mapsto \mbox{tr}_\Omega w \in \mathbb{H}^{\alpha/2}(\Omega),
\end{eqnarray*}
we get the variational formulation of problem \eqref{eq-02}:
to find $\tilde{\mathcal{U}}\in \mathbb{V}$ subject to $\mbox{tr}_\Omega \tilde{\mathcal{U}}|_{t=0}=\psi_0$ and
\begin{eqnarray}\label{eq-03}
  \langle\mbox{tr}_\Omega \partial_t \tilde{\mathcal{U}}, \mbox{tr}_\Omega \phi\rangle_\alpha +
  \frac{1}{d_\alpha}\int_\mathcal{C} (z^\beta \nabla \tilde{\mathcal{U}}) \cdot \nabla \phi
  =\langle f, \mbox{tr}_\Omega \phi\rangle_\alpha, ~~
  \forall \phi \in \mathring{H}^1_L(z^\beta,\mathcal{C}), ~~ \mbox{a.e.}~~ t\in I,
\end{eqnarray}
where $\langle \cdot, \cdot\rangle_\alpha$ is the duality pairing between $\mathbb{H}^{-\alpha/2}(\Omega)$
and $\mathbb{H}^{\alpha/2}(\Omega)$. 

\subsection{The fully implicit discrete FE scheme}

Note that problem \eqref{eq-02} can't be approximated
except that $\mathcal{C}$ is truncated to $\mathcal{C}_\mathcal{Z}:=\Omega\times(0,\mathcal{Z})$
for a suitable $\mathcal{Z}<+\infty$ with
$\partial_L\mathcal{C}_\mathcal{Z}=\partial\Omega\times[0,\mathcal{Z}) \cup \Omega\times\{z=\mathcal{Z}\}$
and an additional homogeneous Dirichlet on $z=\mathcal{Z}$,
whose theoretical foundation is the exponential decay of $\tilde{\mathcal{U}}$ in the extended direction $z$ on the basis of
the length $\mathcal{Z}$, see \cite[Proposition 22]{n-02}.

The first step toward the discretization is to obtain a semi-discrete approximation
of the truncated problem of the original problem \eqref{eq-03}.
To do so, we associate a (possibly nonuniform) mesh by points $0=t_0<t_1<\cdots<t_N=T$
and time-step sizes $\tau_j=t_j-t_{j-1}$ for $j=1,\cdots,N$.
The backward Euler discretization is applied:
to seek $\{U_k\}_{k=0}^N\subset\mathring{H}^1_L(z^\beta,\mathcal{C}_\mathcal{Z})$ subject to $\mbox{tr}_\Omega U_0=\psi_0$ and
\begin{eqnarray*}
  \Big{(}\frac{\mbox{tr}_\Omega U_{k+1}-\mbox{tr}_\Omega U_k}{\tau_{k+1}}, \mbox{tr}_\Omega \vartheta\Big{)}_{L^2(\Omega)} +
  \frac{1}{d_\alpha}\int_{\mathcal{C}_\mathcal{Z}} (z^\beta \nabla U_{k+1}) \cdot \nabla \vartheta
  =\langle f_{k+1}, \mbox{tr}_\Omega \vartheta\rangle_\alpha, ~~k=0,\cdots,N-1
\end{eqnarray*}
for all $\vartheta \in \mathring{H}^1_L(z^\beta,\mathcal{C}_\mathcal{Z})$,
where $f_{k+1} := f|_{t=t_{k+1}}$. Obviously, the sequence $$\{u_k:=\mbox{tr}_\Omega U_k\in \mathbb{H}^{\alpha/2}(\Omega)\}_{k=1}^N$$
is a series of piecewise constant approximations to the solution of problem \eqref{eq-01}.

Next, we turn to discuss the spatial discretization.
Given a quasi-uniform triangulation $\mathcal{J}_\Omega$ over $\Omega$ 
and a graded partition $\mathcal{J}_\mathcal{Z}$ of $[0,\mathcal{Z}]$ by the transition points
\begin{align}
  \label{eq-07}
  z_j=\left \{
    \begin{aligned}
    & z_*\mathcal{Z}\Big{(}\frac{4j}{3M}\Big{)}^\mu, & 0\le j\le \frac{3}{4}M \\
    & \mathcal{Z}\bigg{[}(1-z_*)\Big{(}\frac{4j}{M}-3\Big{)}+z_*\bigg{]}, & \frac{3}{4}M <j\le M
    \end{aligned}
  \right.~~ \mbox{with}~~ \mu=\frac{3}{\alpha}+0.01,~~z_*=\Big{(}1+\frac{\mu}{3}\Big{)}^{-1}
\end{align}
as well as segments $I_k =(z_{k-1},z_k)$ for $k=1,\cdots,M$, we construct the mesh $\mathcal{J}$ over $\mathcal{C}_\mathcal{Z}$
as the tensor product of $\mathcal{J}_\Omega$ and $\mathcal{J}_\mathcal{Z}$, and define
\begin{eqnarray*}
  \mathbb{V}(\mathcal{J}) := \{v\in C(\bar{\mathcal{C}}_\mathcal{Z}):v|_{K\times I_k}\in \mathcal{P}_1(K)\otimes\mathcal{P}_1(I_k),
  \forall K\in \mathcal{J}_\Omega,1\le k\le M; ~v|_{\partial_L\mathcal{C}_\mathcal{Z}} = 0\},~~
  \mathbb{W}(\mathcal{J}_\Omega):=\mbox{tr}_\Omega\mathbb{V}(\mathcal{J}),
\end{eqnarray*}
and the weighted elliptic projector $\mathcal{L}: \mathring{H}^1_L(z^\beta,\mathcal{C}_\mathcal{Z})\rightarrow \mathbb{V}(\mathcal{J})$
such that
\begin{eqnarray*}
  \int_{\mathcal{C}_\mathcal{Z}} \Big{[}z^\beta \nabla (\mathcal{L}w)\Big{]} \cdot \nabla v =
  \int_{\mathcal{C}_\mathcal{Z}} (z^\beta \nabla w) \cdot \nabla v, ~~\forall v\in \mathbb{V}(\mathcal{J})
\end{eqnarray*}
is valid for any function $w\in \mathring{H}^1_L(z^\beta,\mathcal{C}_\mathcal{Z})$,
where $\mathcal{P}_1$ denotes the space of all polynomials of at most one degree.

Finally, we arrive at the fully discrete FE scheme: with the initial guess $U^\mathcal{J}_0:=(\mathcal{L}\circ\mbox{tr}^{-1}_\Omega)\psi_0$,
to calculate $\{U^\mathcal{J}_k\}_{k=1}^N\subset\mathbb{V}(\mathcal{J})$ satisfying the constraint
\begin{eqnarray}\label{eq-04}
  \Big{(}\frac{\mbox{tr}_\Omega U^\mathcal{J}_{k+1}-\mbox{tr}_\Omega U^\mathcal{J}_k}{\tau_{k+1}}, \mbox{tr}_\Omega v\Big{)}_{L^2(\Omega)} +
  \frac{1}{d_\alpha}\int_{\mathcal{C}_\mathcal{Z}} (z^\beta \nabla U^\mathcal{J}_{k+1}) \cdot \nabla v
  =\langle f_{k+1}, \mbox{tr}_\Omega v\rangle_\alpha, ~~\forall v\in \mathbb{V}(\mathcal{J}).
\end{eqnarray}
As before, approximate solutions of problem \eqref{eq-01} are given by
\begin{eqnarray*}
  \{u^\mathcal{J}_k:=\mbox{tr}_\Omega U^\mathcal{J}_k\in\mathbb{W}(\mathcal{J})\}_{k=1}^N.
\end{eqnarray*}

In order to get the matrix formulation of \eqref{eq-04}, it is necessary to introduce 
\begin{eqnarray*}
  \mathbb{U}(\mathcal{J}_\Omega) = \{ u_h\in H_0^1(\Omega): u_h\in C(\Omega),
   u_h|_{K}\in \mathcal{P}_1(K),\forall K\in \mathcal{J}_\Omega \}=\textbf{span}\{\phi_i,~i = 1,\cdots,n\}
\end{eqnarray*}
and
\begin{eqnarray*}
  \mathbb{Z}(\mathcal{J}_\mathcal{Z}) =  \{z_h\in H_0^1(0,\mathcal{Z}): z_h|_{I_k}\in \mathcal{P}_1(I_k),1\le k\le M\}=
  \textbf{span}\{\psi_j,~j=1,\cdots,M-1\},
\end{eqnarray*}
where $n$ is the number of interior vertices in $\mathcal{J}_\Omega$. It is straightforward to verify that
$\mathbb{V}(\mathcal{J})$ and $U^\mathcal{J}_k$ can be respectively written as
\begin{eqnarray*}
  \mathbb{V}(\mathcal{J}) = \textbf{span}\{\phi_i\psi_j,~i = 1,\cdots,n;j =1,\cdots,M-1\}
  ~~ \mbox{and}~~ U^\mathcal{J}_k = \sum_{i=1}^{n} \sum_{j=1}^{M-1} u^{(k)}_{ij}\phi_i\psi_j.
\end{eqnarray*}
Choosing the test function $v = \phi_q\psi_l$ in \eqref{eq-04}, we deduce
\begin{eqnarray*}
  \int_{\mathcal{C}_\mathcal{Z}}\Big{[}(z^\beta\psi_j\nabla\phi_i)\cdot(\phi_q\nabla\psi_l)+
  (z^\beta\phi_i\nabla\psi_j)\cdot(\psi_l\nabla\phi_q)\Big{]}=0
\end{eqnarray*}
by Green's first identity and the fact that $\phi_i|_{\partial\Omega}=\phi_q|_{\partial\Omega}=0$,
and one further expansion via separation of variables
\begin{eqnarray*}
  \frac{1}{d_\alpha}\int_{\mathcal{C}_\mathcal{Z}}\Big{[}z^\beta \nabla (\phi_i\psi_j)\Big{]}\cdot\nabla (\phi_q\psi_l)=a^s_{iq}m^z_{jl}+m^s_{iq}a^z_{jl},
\end{eqnarray*}
where entries
\begin{eqnarray*}
  a^s_{iq} = \frac{1}{d_\alpha}\int_{\Omega}\nabla\phi_i\nabla\phi_q,~~m^z_{jl} = \int_0^{\mathcal Z}z^\beta\psi_j\psi_l,
  ~~m^s_{iq} = \frac{1}{d_\alpha}\int_{\Omega}\phi_i\phi_q~~\mbox{and}~~a^z_{jl} = \int_0^{\mathcal Z} z^\beta\frac{d\psi_j}{d z}\frac{d\psi_l}{dz}.
\end{eqnarray*}
Thus, the system of linear equations obtained from \eqref{eq-04} is of the form
\begin{equation}\label{eq-08}
\left(\frac{1}{\tau_{k+1}}\begin{bmatrix}
  M_s & O \\
  O & O
 \end{bmatrix}
 + M_z\otimes A_s+A_z\otimes M_s\right)\begin{bmatrix}
\mathcal{U}_{k+1} \\
\mathcal{V}_{k+1}
\end{bmatrix}
=\frac{1}{\tau_{k+1}}\begin{bmatrix}
  M_s & O \\
  O & O
 \end{bmatrix}\begin{bmatrix}
\mathcal{U}_k\\
\mathcal{V}_k
\end{bmatrix}+\begin{bmatrix}
\mathcal{F}_{k+1}\\
O
\end{bmatrix},~~k=0,\cdots,N-1,
\end{equation}
where
\begin{itemize}

  \item matrixes $M_z = (m^z_{ij})_{(M-1)\times(M-1)}$, $A_z = (a^z_{ij})_{(M-1)\times(M-1)}$,
  $M_s = (m^s_{ij})_{n\times n}$ and $A_s = (a^s_{ij})_{n\times n}$,

  \item $O$ represents a zero matrix of suitable size, $\otimes$ denotes the Kronecker product,

  \item $\mathcal{U}_k=(u^{(k)}_{11},\cdots,u^{(k)}_{1n})^T$ with $u^{(0)}_{1i}=\tilde{\psi}_0(x_i,y_i)$ for $i=1,\cdots,n$
to approximate the initial $\psi_0(x,y)$,

  \item $\mathcal{V}_k=(u^{(k)}_{21},\cdots,u^{(k)}_{2n};\cdots;u^{(k)}_{(M-1)1},\cdots,u^{(k)}_{(M-1)n})^T$
is the auxiliary vector produced by the other points in the extended direction,

  \item $\mathcal{F}_{k+1}=(f^{(k+1)}_1,\cdots,f^{(k+1)}_n)^T$ with
$f^{(k+1)}_i=\langle f_{k+1}, \phi_i\rangle_\alpha$.

\end{itemize}
It is trivial to see that the coefficient matrix of \eqref{eq-08} is symmetric positive definite (SPD).
Therefore, the multigrid with line smoother proposed in \cite{c-01} is a very suitable method and is employed in our implementation.

\section{MGRIT formulation and a generalized two-level convergence upper bound}\label{sec3}

The main focus of this section is the general MGRIT algorithm using FCF-relaxation,
followed by its two-level convergence analysis for the general case
where the time-grid propagators don't have to be the same and need not be simultaneously
diagonalizable with a unitary matrix.

\subsection{The general MGRIT algorithm using FCF-relaxation} \label{sec3-1}

Consider the following block unit lower bidiagonal system
\begin{equation}\label{eq-05}
\mathcal{A}\mathcal{U}:=
\begin{bmatrix}
  I &        &        &   \\
  -\Psi_1 & I &        &   \\
        & \ddots & \ddots &   \\
        &        &  -\Psi_N & I
\end{bmatrix}
\begin{bmatrix}
\mathcal{U}_0 \\
\mathcal{U}_1 \\
\vdots \\
\mathcal{U}_N
\end{bmatrix}
=
\begin{bmatrix}
\mathcal{G}_0 \\
\mathcal{G}_1 \\
\vdots \\
\mathcal{G}_N
\end{bmatrix}:=\mathcal{G},
\end{equation}
where $\mathcal{G}_0=\mathcal{U}_0$, $I$ is the identity matrix,
vectors $\{\mathcal{U}_j\}_{j=1}^N$, $\{\mathcal{G}_j\}_{j=1}^N$ and time-grid propagators $\{\Psi_j\}_{j=1}^N$
are assumed not to be connected with any specific forms.

\begin{figure}[htbp]
\centerline{
\includegraphics[scale=0.3]{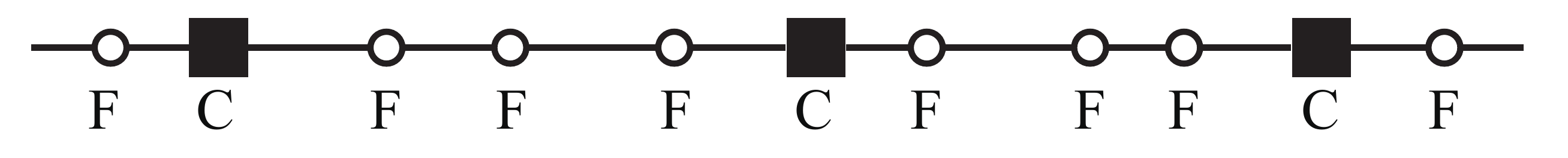}}
\caption{Schematic view of the CF-splitting by a factor of four.}\label{fig-01}
\end{figure}

Various components are required in order to solve the global space-time problem \eqref{eq-05} through the MGRIT algorithm
based on the multigrid reduction technique \cite{r-01}.
Define the coarse temporal mesh with nodes $\{\tilde{t}_i:=t_{m\times i}\}_{i=0}^{N_c}$
and time-step sizes $\{\tilde{\tau}_j=\tilde{t}_j-\tilde{t}_{j-1}\}_{j=1}^{N_c}$ for some coarsening factor $m$,
where $N_c=N/m$. Then all $\tilde{t}_i$ are disconnected C-points representing those variables in the coarse level
and the others are F-points grouped $m-1$ adjacent points into a bunch.
The particular case when $m=4$ is illustrated in Fig. \ref{fig-01}.
It should be mentioned that $\{\tilde{t}_i\}_{i=0}^{N_c}$ is still viewed to be nonuniform.
Assume for simplicity of presentation here that F-points are ordered first.
MGRIT utilizes the coarse injection $[0,I_c]$ as the restriction operator,
i.e., only the residual at C-points is computed, where $I_c$ is the identity matrix.
The FCF-relaxation (i.e., an initial F-relaxation followed by a C-relaxation and then a second F-relaxation)
is often the most preferable choice to produce optimal and scalable multilevel iterations \cite{f-02},
where F-relaxation refers to updating variables at $m-1$ F-points in turn yet independent of other F-intervals,
and C-relaxation stands for performing \eqref{eq-08} to propagate the solution at a C-point. It is apparent that these two processes are both highly parallel.
The ideal interpolation is used to provide the coarse-grid correction,
which is implemented at the cost of updating C-variables by the coarse approximations immediately
followed by an F-relaxation over the time-grid propagators $\{\Psi_j\}_{j=1}^N$.
It is worth mentioning that the resulting Petrov-Galerkin triple matrix product is exactly the Schur complement
\begin{equation}\label{eq-10}
\mathcal{A}_\Delta=\begin{bmatrix}
I               &                 &          &                   &   \\
-\mathcal{B}_{\Psi,m}^m & I               &          &                   &   \\
                & -\mathcal{B}_{\Psi,2m}^m & I        &                   &   \\
                &                 & \ddots   & \ddots            &    \\
                &                 &          & -\mathcal{B}_{\Psi,N}^m & I
\end{bmatrix},
\end{equation}
where $\mathcal{B}_{\Psi,s}^l=\Psi_s\Psi_{s-1}\cdots\Psi_{s-l+1}$ with $\mathcal{B}_{\Psi,s}^0=I$,
positive integer $s$ and nonnegative integer $l$. However,
this coarse-grid problem would be as expensive as the original fine-grid problem.
MGRIT approximates the coarse-grid operator with
\begin{equation}\label{eq-11}
\mathcal{S}_\Delta=\begin{bmatrix}
I               &                 &          &                   &   \\
-\Phi_1 & I               &          &                   &   \\
                & -\Phi_2 & I        &                   &   \\
                &                 & \ddots   & \ddots            &    \\
                &                 &          & -\Phi_{N_c} & I
\end{bmatrix},
\end{equation}
where $\Phi_k$ $(k=1,\cdots,N_c)$ is an approximate coarse-grid time-stepping operator to $\mathcal{B}_{\Psi,km}^m$.
The most natural choice of $\Phi_k$ is to re-discretize the time-dependent linear continuous problem 
on the coarse temporal mesh.
The final stage is the application of the above processes recursively, with possibly different coarsening factors at different time levels,
to construct a multilevel hierarchy under certain stopping criteria, e.g.,
by setting the maximum number of time levels and the minimum coarsest temporal grid size.
Below is a sketch of the resulting MGRIT V(1,0)-cycle algorithm using FCF-relaxation in the recursive fashion,
where $L$ is the number of time levels determined by the user supplied stopping procedure,
$\mathcal{A}^{(1)}=\mathcal{A}$, $\mathcal{G}^{(1)}=\mathcal{G}$,
matrices $\mathcal{A}^{(k+1)}$, $\mathcal{R}^{(k)}$ and $\mathcal{P}^{(k)}$
($k=1,\cdots,L-1$) respectively correspond to the $k$-th coarse-scale time re-discretization,
the mentioned above restriction and interpolation.

\begin{Algorithm}\label{alg-02}
MGRIT V(1,0)-cycle algorithm using FCF-relaxation:
$\mathcal{U}^{(k)}=\textrm{MGRIT}(k,\mathcal{A}^{(k)},\mathcal{U}^{(k)},\mathcal{G}^{(k)})$
\begin{description}

  \item[Step 1] Apply FCF-relaxation to $\mathcal{A}^{(k)} \mathcal{U}^{(k)} = \mathcal{G}^{(k)}$ with $\mathcal{U}^{(k)}$ as the initial guess.

  \item[Step 2] Inject the residual to the coarse time level
  $\mathcal{G}^{(k+1)}=\mathcal{R}^{(k)}(\mathcal{G}^{(k)} - \mathcal{A}^{(k)} \mathcal{U}^{(k)})$.

  \item[Step 3] If $k+1=L$, then solve $\mathcal{A}^{(L)}\mathcal{U}^{(L)}=\mathcal{G}^{(L)}$; Else compute
  $\mathcal{U}^{(k+1)}=\textrm{MGRIT}(k+1,\mathcal{A}^{(k+1)},0,\mathcal{G}^{(k+1)})$.

  \item[Step 4] Perform the coarse-grid correction using the ideal interpolation
  $\mathcal{U}^{(k)} = \mathcal{U}^{(k)} + \mathcal{P}^{(k)} \mathcal{U}^{(k+1)}$.

\end{description}
\end{Algorithm}

\begin{remark} (Reduction of the computational overhead)
  Here, we point out that
  Step 4 of Algorithm \ref{alg-02} is executed by just correcting C-variables,
  and updating F-variables only when the Euclidean norm of the residual is small enough.
  The underlying reason is that corrections at F-points are equivalent to an F-relaxation,
  which will be performed in the subsequent MGRIT iteration.
\end{remark}

\subsection{A generalized two-level convergence analysis}\label{sec3-2}

Setting $L=2$, Algorithm \ref{alg-02} in this case reduces to the well-known two-level scheme,
whose residual propagation operator for a single iteration, as in the derivations of \cite{y-01,s-03},
can be expressed as in an FC-partitioning of $\mathcal{A}$:
\begin{eqnarray}\label{eq-14}
 \mathcal{R}_{TL}^{FCF}=\begin{bmatrix}
  0 \\
  \mathcal{T}_\Psi^\Phi
  \end{bmatrix}[\mathcal{Q} ~~ I],\quad
  \mathcal{T}_\Psi^\Phi=(I-\mathcal{A}_\Delta \mathcal{S}_\Delta^{-1})(I-\mathcal{A}_\Delta),
\end{eqnarray}
where the matrix used in forming the ``ideal restriction" operator
\begin{eqnarray*}
\mathcal{Q}=\left[
  \begin{array}{ccccccccccccc}
  0 & \cdots & 0 & 0 &  &  &  &  &  &  &  &  &  \\
  \mathcal{B}_{\Psi,m}^{m-1} & \cdots & \mathcal{B}_{\Psi,m}^2 & \mathcal{B}_{\Psi,m}^1 &  &  &  &  &  &  &  &  &  \\
  &  &  &  & \mathcal{B}_{\Psi,2m}^{m-1} & \cdots & \mathcal{B}_{\Psi,2m}^2 & \mathcal{B}_{\Psi,2m}^1 &  &  &  &  &  \\
  &  &  &  &  &  &  &  & \ddots &  &  &  &  \\
  &  &  &  &  &  &  &  &  & \mathcal{B}_{\Psi,N}^{m-1} & \cdots & \mathcal{B}_{\Psi,N}^2 & \mathcal{B}_{\Psi,N}^1
\end{array}\right].
\end{eqnarray*}
Below, we wish to find an upper bound for a certain norm of the residual propagation $\mathcal{R}_{TL}^{FCF}$.
\begin{description}

  \item[S1] Assume that time-stepping operators $\Psi_j$ and $\Phi_k$ are both strongly stable \cite{l-04}
(i.e., the discrete $l^2$-norms $\|\Psi_j\|<1$ and $\|\Phi_k\|<1$) and simultaneously diagonalizable by the matrix $\mathcal{X}_s$.
Note that the assumption on simultaneous diagonalization will be tenable
when the same spatial discretization is exploited on fine and coarse temporal meshes.
Denote by $\{\lambda^{(\omega)}_j\}_{\omega=1}^n$ and $\{\mu^{(\omega)}_k\}_{\omega=1}^n$
the sets of all complex eigenvalues of $\Psi_j$ and $\Phi_k$, respectively.

  \item[S2] Through simple algebraic calculations, we arrive at an estimate of the second term in \eqref{eq-14}:
\begin{eqnarray*}
 \big{\|}[\mathcal{Q} ~~ I]\big{\|}=\sqrt{\lambda_{\textrm{max}}\left(\begin{bmatrix}
  \mathcal{Q}^T \\
  I
  \end{bmatrix}[\mathcal{Q} ~~ I]\right)}
 =\max_k\sqrt{\left\|\sum_{i=0}^{m-1}\mathcal{B}_{\Psi,km}^i(\mathcal{B}_{\Psi,km}^i)^T\right\|}\le
 \max_k\sqrt{\sum_{i=0}^{m-1}\left\|\mathcal{B}_{\Psi,km}^i\right\|^2}.
\end{eqnarray*}
From the relation $\big{\|}\Psi_j\big{\|}<1$ for $j=1,\cdots,N$ in Stage {\bf S1}, we have
\begin{eqnarray*}
 \left\|\mathcal{B}_{\Psi,km}^0\right\|=1,~~
 \left\|\mathcal{B}_{\Psi,km}^i\right\|\le\big{\|}\Psi_{km}\big{\|}\cdot
 \big{\|}\Psi_{km-1}\big{\|}\cdots\big{\|}\Psi_{km-i+1}\big{\|}<1~\mbox{for}~i=1,\cdots,m-1,
\end{eqnarray*}
and hence $\big{\|}[\mathcal{Q} ~~ I]\big{\|}<\sqrt{m}$.

  \item[S3] By exploiting the closed form
\begin{equation*}
\mathcal{S}_\Delta^{-1}=\begin{bmatrix}
I               &                 &          &                   &   \\
\mathcal{B}_{\Phi,1}^1 & I               &          &                   &   \\
\mathcal{B}_{\Phi,2}^2 & \mathcal{B}_{\Phi,2}^1 & I        &                   &   \\
\vdots & \vdots & \ddots   & \ddots            &    \\
\mathcal{B}_{\Phi,N_c}^{N_c} & \mathcal{B}_{\Phi,N_c}^{N_c-1} & \cdots & \mathcal{B}_{\Phi,N_c}^1 & I
\end{bmatrix},
\end{equation*}
the other two matrix forms appeared in $\mathcal{R}_{TL}^{FCF}$ are
obtained from block lower bidiagonal structures \eqref{eq-10}-\eqref{eq-11}:
\begin{equation*}
I-\mathcal{A}_\Delta \mathcal{S}_\Delta^{-1}=\begin{bmatrix}
0               &                 &          &                   &   \\
\mathcal{B}_{\Psi,m}^m-\Phi_1 & 0               &          &                   &   \\
(\mathcal{B}_{\Psi,2m}^m-\Phi_2)\mathcal{B}_{\Phi,1}^1 & \mathcal{B}_{\Psi,2m}^m-\Phi_2 & 0        &                   &   \\
\vdots & \vdots & \ddots   & \ddots            &    \\
(\mathcal{B}_{\Psi,N}^m-\Phi_{N_c})\mathcal{B}_{\Phi,N_c-1}^{N_c-1}
& (\mathcal{B}_{\Psi,N}^m-\Phi_{N_c})\mathcal{B}_{\Phi,N_c-1}^{N_c-2} & \cdots & \mathcal{B}_{\Psi,N}^m-\Phi_{N_c} & 0
\end{bmatrix}
\end{equation*}
and
\begin{equation*}
I-\mathcal{A}_\Delta=\begin{bmatrix}
0               &                 &          &                   &   \\
\mathcal{B}_{\Psi,m}^m & 0               &          &                   &   \\
                & \mathcal{B}_{\Psi,2m}^m & 0        &                   &   \\
                &                 & \ddots   & \ddots            &    \\
                &                 &          & \mathcal{B}_{\Psi,N}^m & 0
\end{bmatrix}.
\end{equation*}
As a result, the expression of $\mathcal{T}_\Psi^\Phi$ can be rewritten as
\begin{equation}\label{eq-15}
\mathcal{T}_\Psi^\Phi=\begin{bmatrix}
0               &                 &          &                   &  & \\
0               &      0          &          &                   &  & \\
(\mathcal{B}_{\Psi,2m}^m-\Phi_2)\mathcal{B}_{\Psi,m}^m & 0               &    0    &                   &  & \\
(\mathcal{B}_{\Psi,3m}^m-\Phi_3)\mathcal{B}_{\Phi,2}^1\mathcal{B}_{\Psi,m}^m
& (\mathcal{B}_{\Psi,3m}^m-\Phi_3)\mathcal{B}_{\Psi,2m}^m &    0      &      0           &  &  \\
  \vdots              &    \vdots     &   \ddots  &     \ddots       & \ddots  & \\
(\mathcal{B}_{\Psi,N}^m-\Phi_{N_c})\mathcal{B}_{\Phi,N_c-1}^{N_c-2}\mathcal{B}_{\Psi,m}^m
& (\mathcal{B}_{\Psi,N}^m-\Phi_{N_c})\mathcal{B}_{\Phi,N_c-1}^{N_c-3}\mathcal{B}_{\Psi,2m}^m & \cdots
& (\mathcal{B}_{\Psi,N}^m-\Phi_{N_c})\mathcal{B}_{\Psi,N-m}^m & 0 & 0
\end{bmatrix},
\end{equation}
which can be reformulated as the under-mentioned triple matrix product
\begin{equation*}
\begin{bmatrix}
I &  &   &  &  &  \\
  & I &  &   &  &  \\
  &  & \mathcal{B}_{\Psi,2m}^m-\Phi_2 &   &  &  \\
&  &  & & \ddots &  \\
&  &  & &  & \mathcal{B}_{\Psi,N}^m-\Phi_{N_c}
\end{bmatrix}
\begin{bmatrix}
0 &   &   &   &  &  \\
0 & 0 &   &   &  &  \\
I & 0 & 0 &   &  &  \\
\mathcal{B}_{\Phi,2}^1 & I & 0 & 0 &  &  \\
\vdots & \vdots & \ddots & \ddots & \ddots &  \\
\mathcal{B}_{\Phi,N_c-1}^{N_c-2} & \mathcal{B}_{\Phi,N_c-1}^{N_c-3} & \cdots & I & 0 & 0
\end{bmatrix}
\begin{bmatrix}
\mathcal{B}_{\Psi,m}^m &  &   & & \\
& \ddots &  & & \\
&  &  \mathcal{B}_{\Psi,N-m}^m & & \\
&  &  & I & \\
&  &  &  & I
\end{bmatrix}.
\end{equation*}
Then, under the assumption that $\mathcal{B}_{\Psi,jm}^m-\Phi_j$ is invertible for $j=2,\cdots,N_c$,
using the relation $\mathcal{B}_{\Phi,s}^k=\mathcal{B}_{\Phi,s}^l\mathcal{B}_{\Phi,s-l}^{k-l}$,
and recalling that the Moore-Penrose conditions are used to define the pseudoinverse $D^\dag$ of any matrix $D$:
\begin{equation*}
 DD^\dag~\mbox{and}~D^\dag D~\mbox{are both Hermitian}, DD^\dag D=D~\mbox{and}~D^\dag D D^\dag = D^\dag,
\end{equation*}
the unique pseudoinverse $(\mathcal{T}_\Psi^\Phi)^\dag$ of $\mathcal{T}_\Psi^\Phi$
can be established in a straightforward way
\begin{equation}\label{eq-16}
(\mathcal{T}_\Psi^\Phi)^\dag=\begin{bmatrix}
0 & 0 & t_2^\dag &          &                &    \\
0 & 0 & g_2^\dag & t_3^\dag &                &    \\
  &   &          & \ddots   & \ddots         &    \\
  &   &          &          & g_{N_c-1}^\dag & t_{N_c}^\dag \\
  &   &          &          & 0              & 0  \\
  &   &          &          & 0              & 0
\end{bmatrix}~~\mbox{with}~~
\begin{matrix}
   t_l^\dag = \big{(}\mathcal{B}_{\Psi,(l-1)m}^m\big{)}^{-1}(\mathcal{B}_{\Psi,lm}^m-\Phi_l)^{-1},~l=2,\cdots,N_c, \\
   \\
   g_s^\dag = -(\mathcal{B}_{\Psi,sm}^m)^{-1}\Phi_s(\mathcal{B}_{\Psi,sm}^m-\Phi_s)^{-1},~s=2,\cdots,N_c-1,
\end{matrix}
\end{equation}
by noticing the following two important results
\begin{equation*}
\mathcal{T}_\Psi^\Phi(\mathcal{T}_\Psi^\Phi)^\dag=\begin{bmatrix}
0 &   &   &  &  \\
 & 0 &   &  &  \\
 &  & I &  &  \\
 &  &  & \ddots &  \\
 &  &  &  & I
\end{bmatrix},\quad
(\mathcal{T}_\Psi^\Phi)^\dag\mathcal{T}_\Psi^\Phi=\begin{bmatrix}
I &  &   & & \\
& \ddots &  & & \\
&  &  I & & \\
&  &  & 0 & \\
&  &  &  & 0
\end{bmatrix}.
\end{equation*}
It is easy to see that the use of $(\mathcal{T}_\Psi^\Phi)^\dag$
allows for a simpler technical path to analyze the maximum singular value of $\mathcal{T}_\Psi^\Phi$
through the minimum nonzero singular value of $(\mathcal{T}_\Psi^\Phi)^\dag$ with the advantage of its simplicity.

  \item[S4] Based on the simultaneous diagonalizability from Stage {\bf S1}: 
\begin{eqnarray*}
 \mathcal{X}_s^{-1}\Psi_j\mathcal{X}_s=\textbf{diag}(\lambda^{(1)}_j,\cdots,\lambda^{(n)}_j):=\Lambda_j,~j=1,\cdots,N;~~
 \mathcal{X}_s^{-1}\Phi_k\mathcal{X}_s=\textbf{diag}(\mu^{(1)}_k,\cdots,\mu^{(n)}_k):=\Upsilon_k,~k=1,\cdots,N_c,
\end{eqnarray*}
we have
\begin{eqnarray*}
\big{\|}\mathcal{T}_\Psi^\Phi\big{\|}_{(\tilde{\mathcal{X}}_s\tilde{\mathcal{X}}_s^*)^{-1}}=
\sup_{v\ne0}\frac{\big{\|}\tilde{\mathcal{X}}_s^{-1}\mathcal{T}_\Psi^\Phi\tilde{\mathcal{X}}_s
\tilde{\mathcal{X}}_s^{-1}v\big{\|}}{\big{\|}\tilde{\mathcal{X}}_s^{-1}v\big{\|}}=
\sup_{w\ne0}\frac{\big{\|}\mathcal{P}_s\mathcal{T}_\Lambda^\Upsilon\mathcal{P}_s^T
w\big{\|}}{\big{\|}w\big{\|}}=
\max_\omega\big{\|}\big{[}\mathcal{T}_\Lambda^\Upsilon\big{]}_\omega\big{\|}=
\frac{1}{\min\limits_\omega \sigma_{\min}\big{(}\big{[}(\mathcal{T}_\Lambda^\Upsilon)^\dag\big{]}_\omega\big{)}},
\end{eqnarray*}
where
\begin{itemize}

  \item $\tilde{\mathcal{X}}_s=\textbf{diag}(\mathcal{X}_s,\cdots,\mathcal{X}_s)$ consisting of $N_c+1$ diagonal blocks,

  \item $\mathcal{T}_\Lambda^\Upsilon$ and $(\mathcal{T}_\Lambda^\Upsilon)^\dag$
  are obtained from replacing $\Psi$, $\Phi$ by $\Lambda$, $\Upsilon$ in \eqref{eq-15} and \eqref{eq-16}, respectively,

  \item $\mathcal{P}_s$ denotes the orthogonal permutation matrix
  which makes $\mathcal{P}_s\mathcal{T}_\Lambda^\Upsilon\mathcal{P}_s^T$ block diagonal,

  \item $\big{[}\mathcal{T}_\Lambda^\Upsilon\big{]}_\omega$ and
  $\big{[}(\mathcal{T}_\Lambda^\Upsilon)^\dag\big{]}_\omega$ represent $\mathcal{T}_\Lambda^\Upsilon$ and
  $(\mathcal{T}_\Lambda^\Upsilon)^\dag$, respectively, evaluated at $\lambda^{(\omega)}_j$ ($j=1,\cdots,N$) and $\mu^{(\omega)}_k$ ($k=1,\cdots,N_c$).

\end{itemize}
We need to determine the nonzero singular value $\sigma_{\min}\big{(}\big{[}(\mathcal{T}_\Lambda^\Upsilon)^\dag\big{]}_\omega\big{)}$.
For this purpose, deleting the last two zero-rows and zero-columns of
$\big{[}(\mathcal{T}_\Lambda^\Upsilon)^\dag\big{]}_\omega\big{[}(\mathcal{T}_\Lambda^\Upsilon)^\dag\big{]}_\omega^*$
which correspond to two zero eigenvalues, the remaining algebraic manipulation is to pursue after
the square root of the smallest nonzero eigenvalue of the complex Hermitian tridiagonal matrix
\begin{eqnarray}\label{eq-18}
\begin{bmatrix}
 \zeta_2^{(\omega)} & -\rho_2^{(\omega)} &  &  &  &  \\
 -\bar{\rho}_2^{(\omega)} & \theta_2^{(\omega)}+\zeta_3^{(\omega)} & -\rho_3^{(\omega)} &  &  &  \\
         & \ddots & \ddots & \ddots &  &  \\
 &  &  & -\bar{\rho}_{N_c-2}^{(\omega)} & \theta_{N_c-2}^{(\omega)}+\zeta_{N_c-1}^{(\omega)} & -\rho_{N_c-1}^{(\omega)} \\
 &  &  &  & -\bar{\rho}_{N_c-1}^{(\omega)} & \theta_{N_c-1}^{(\omega)}+\zeta_{N_c}^{(\omega)}
\end{bmatrix}
\end{eqnarray}
with
\begin{align*}
    \begin{aligned}
    & \zeta_\epsilon^{(\omega)}=\frac{1}{\left|\lambda_{(\epsilon-1)m}^{(\omega)}\cdots\lambda_{(\epsilon-2)m+1}^{(\omega)}\right|^2
 \left|\lambda_{\epsilon m}^{(\omega)}\cdots\lambda_{(\epsilon-1)m+1}^{(\omega)}-\mu_\epsilon^{(\omega)}\right|^2},\quad\epsilon=2,\cdots,N_c; \\
    & \rho_\varsigma^{(\omega)}=\frac{\bar{\mu}_\varsigma^{(\omega)}}{\bar{\lambda}_{\varsigma m}^{(\omega)}
 \cdots\bar{\lambda}_{(\varsigma-1)m+1}^{(\omega)}\lambda_{(\varsigma-1)m}^{(\omega)}\cdots\lambda_{(\varsigma-2)m+1}^{(\omega)}
 \left|\lambda_{\varsigma m}^{(\omega)}\cdots\lambda_{(\varsigma-1)m+1}^{(\omega)}-\mu_\varsigma^{(\omega)}\right|^2}, \\
    & \theta_\varsigma^{(\omega)}=\frac{\left|\mu_\varsigma^{(\omega)}\right|^2}{\left|\lambda_{\varsigma m}^{(\omega)}
    \cdots\lambda_{(\varsigma-1)m+1}^{(\omega)}\right|^2
 \left|\lambda_{\varsigma m}^{(\omega)}\cdots\lambda_{(\varsigma-1)m+1}^{(\omega)}-\mu_\varsigma^{(\omega)}\right|^2},\qquad\varsigma=2,\cdots,N_c-1.
    \end{aligned}
\end{align*}
It is assumed throughout that the temporal eigenvalue approximation property (TEAP) is valid for each eigenmode $\omega$
with positive constants $\{\delta_\varsigma^{(\omega)}\}_{\varsigma=2}^{N_c-1}$:
\begin{eqnarray}\label{eq-20}
 \left|\lambda_{\varsigma m}^{(\omega)}\cdots\lambda_{(\varsigma-1)m+1}^{(\omega)}-\mu_\varsigma^{(\omega)}\right|^2\le
 \delta_\varsigma^{(\omega)}\Big{(}\frac{1}{\left|\lambda_{(\varsigma-1)m}^{(\omega)}\cdots\lambda_{(\varsigma-2)m+1}^{(\omega)}\right|^2}
 -\frac{\left|\mu_\varsigma^{(\omega)}\right|}{\left|\lambda_{\varsigma m}^{(\omega)}\cdots\lambda_{(\varsigma-2)m+1}^{(\omega)}\right|}\Big{)},
 ~~\varsigma=2,\cdots,N_c-1.
\end{eqnarray}
Applying the Gershgorin circle theorem, it can be deduced that
\begin{eqnarray}\label{eq-17}
\big{\|}\mathcal{T}_\Psi^\Phi\big{\|}^2_{(\tilde{\mathcal{X}}_s\tilde{\mathcal{X}}_s^*)^{-1}}\le\max_\omega\max\left\{\delta_2^{(\omega)},
\max_{\varsigma=2,\cdots,N_c-1}\frac{\delta_{\varsigma+1}^{(\omega)}\delta_\varsigma^{(\omega)}
\left|\lambda_{\varsigma m}^{(\omega)}\cdots\lambda_{(\varsigma-1)m+1}^{(\omega)}\right|}{\delta_\varsigma^{(\omega)}
\left|\lambda_{\varsigma m}^{(\omega)}\cdots\lambda_{(\varsigma-1)m+1}^{(\omega)}\right|-
\delta_{\varsigma+1}^{(\omega)}\left|\mu_\varsigma^{(\omega)}\right|
\left|\lambda_{(\varsigma-1)m}^{(\omega)}\cdots\lambda_{(\varsigma-2)m+1}^{(\omega)}\right|}\right\}
\end{eqnarray}
by making the crucial observations that
\begin{eqnarray*}
 \frac{1}{\left|\lambda_{\varsigma m}^{(\omega)}\cdots\lambda_{(\varsigma-1)m+1}^{(\omega)}\right|^2
 \left|\lambda_{(\varsigma+1)m}^{(\omega)}\cdots\lambda_{\varsigma m+1}^{(\omega)}-\mu_{\varsigma+1}^{(\omega)}\right|^2}-
 \frac{\left|\mu_{\varsigma+1}^{(\omega)}\right|}{\left|\lambda_{(\varsigma+1)m}^{(\omega)}\cdots\lambda_{(\varsigma-1)m+1}^{(\omega)}\right|
 \left|\lambda_{(\varsigma+1)m}^{(\omega)}\cdots\lambda_{\varsigma m+1}^{(\omega)}-\mu_{\varsigma+1}^{(\omega)}\right|^2}
 \ge\frac{1}{\delta_\varsigma^{(\omega)}}
\end{eqnarray*}
and
\begin{eqnarray*}
\frac{1}{\left|\lambda_{\varsigma m}^{(\omega)}\cdots\lambda_{(\varsigma-1)m+1}^{(\omega)}-\mu_\varsigma^{(\omega)}\right|^2}\Big{(}
\frac{\left|\mu_\varsigma^{(\omega)}\right|^2}{\left|\lambda_{\varsigma m}^{(\omega)}\cdots\lambda_{(\varsigma-1)m+1}^{(\omega)}\right|^2}-
\frac{\left|\mu_\varsigma^{(\omega)}\right|}{\left|\lambda_{\varsigma m}^{(\omega)}\cdots\lambda_{(\varsigma-2)m+1}^{(\omega)}\right|}\Big{)}
\ge-\frac{\left|\mu_\varsigma^{(\omega)}\right|}{\delta_\varsigma^{(\omega)}}\times
\frac{\left|\lambda_{(\varsigma-1)m}^{(\omega)}\cdots\lambda_{(\varsigma-2)m+1}^{(\omega)}\right|}{\left|\lambda_{\varsigma m}^{(\omega)}
\cdots\lambda_{(\varsigma-1)m+1}^{(\omega)}\right|}.
\end{eqnarray*}
\end{description}
These four stages are summarized in the following theorem.

\begin{theorem}\label{thm-01}
 Assume that $\Psi_j$ and $\Phi_k$ are both strongly stable
 and simultaneously diagonalizable by $\mathcal{X}_s$ with eigenvalues
 $\{\lambda^{(\omega)}_j\}_{\omega=1}^n$ and $\{\mu^{(\omega)}_k\}_{\omega=1}^n$, respectively.
 Let $\sigma_{\min}^{(\omega)}$ be the smallest nonzero eigenvalue of \eqref{eq-18}. Then,
\begin{eqnarray}\label{eq-21}
\big{\|}\mathcal{T}_\Psi^\Phi\big{\|}_{(\tilde{\mathcal{X}}_s\tilde{\mathcal{X}}_s^*)^{-1}}=
\frac{1}{\min\limits_\omega\sqrt{\sigma_{\min}^{(\omega)}}}.
\end{eqnarray}
 If the TEAP condition \eqref{eq-20} holds for all $\varsigma=2,\cdots,N_c-1$ and $\omega=1,\cdots,n$,
 then \eqref{eq-17} provides an actual upper bound. Furthermore,
\begin{eqnarray*}
\delta_2^{(\omega)} < 1 ~~\mbox{and}~~
\frac{\delta_{\varsigma+1}^{(\omega)}\delta_\varsigma^{(\omega)}
\left|\lambda_{\varsigma m}^{(\omega)}\cdots\lambda_{(\varsigma-1)m+1}^{(\omega)}\right|}{\delta_\varsigma^{(\omega)}
\left|\lambda_{\varsigma m}^{(\omega)}\cdots\lambda_{(\varsigma-1)m+1}^{(\omega)}\right|-
\delta_{\varsigma+1}^{(\omega)}\left|\mu_\varsigma^{(\omega)}\right|
\left|\lambda_{(\varsigma-1)m}^{(\omega)}\cdots\lambda_{(\varsigma-2)m+1}^{(\omega)}\right|}<1
\end{eqnarray*}
 are the sufficient conditions to ensure the convergence of the two-level version of Algorithm 1.
\end{theorem}

\begin{remark}
Notice that, for the particular case of the equidistant fine temporal mesh,
\eqref{eq-18} can be written as
\begin{eqnarray*}
\frac{1}{\left|(\lambda^{(\omega)})^m\right|^2
 \left|(\lambda^{(\omega)})^m-\mu^{(\omega)}\right|^2}\mathcal{D}_{\omega},~~
 \mbox{where}~~\mathcal{D}_{\omega}=\begin{bmatrix}
 1 & -\bar{\mu}^{(\omega)} &  &  &  &  \\
 -\mu^{(\omega)} & 1+|\mu^{(\omega)}|^2 & -\bar{\mu}^{(\omega)} &  &  &  \\
         & \ddots & \ddots & \ddots &  &  \\
 &  &  & -\mu^{(\omega)} & 1+|\mu^{(\omega)}|^2 & -\bar{\mu}^{(\omega)} \\
 &  &  &  & -\mu^{(\omega)} & 1+|\mu^{(\omega)}|^2
\end{bmatrix}.
\end{eqnarray*}
Following from the general result \cite[equality (7)]{y-03} and the necessary condition \cite[equality (6)]{y-03},
it is not difficult to see that
the smallest eigenvalue of $\mathcal{D}_{\omega}$, denoted by $\lambda_{\min}(\mathcal{D}_{\omega})$, is bounded by
\begin{eqnarray*}
(1-|\mu^{(\omega)}|)^2+\frac{\pi^2|\mu^{(\omega)}|}{6(N_c-1)^2}\le\lambda_{\min}(\mathcal{D}_{\omega})
\le(1-|\mu^{(\omega)}|)^2+\frac{\pi^2|\mu^{(\omega)}|}{(N_c-1)^2},
\end{eqnarray*}
which can be derived from an argument analogous to that for \cite[Lemma 36]{s-03},
where perturbations $\tilde{\alpha}=0$ and $\tilde{\beta}=|\mu^{(\omega)}|^2$
are replaced by $\tilde{\alpha}=|\mu^{(\omega)}|^2$ and $\tilde{\beta}=0$. Indeed, in this case,
the upper bound \eqref{eq-21} would be reduced to the tight bounds
stated in \cite[Theorem 30, inequality (63)]{s-03}.
Furthermore, it is apparent that Theorem \ref{thm-01} can be regarded as a generalization of \cite[Theorem 35]{s-03}
to the FCF-relaxation case, which has been proven a key component used to build a scalable multilevel solver \cite{f-02}.
\end{remark}

\begin{remark}
The $(\tilde{\mathcal{X}}_s\tilde{\mathcal{X}}_s^*)^{-1}$-norm will be reduced to
the discrete $l^2$-norm only when the matrix $\mathcal{X}_s$ is unitary, that is,
$\Psi_j$ and $\Phi_k$ are both normal.
\end{remark}

\begin{remark}
The upper bound \eqref{eq-17} is usually not tight. Fortunately,
the exact convergence upper bound \eqref{eq-21} can be evaluated
via hunting for the smallest eigenvalue of \eqref{eq-18}.
Multitudinous research efforts on open-source packages
have been made to deal with this issue, e.g., LAPACK, SLEPc, EISPACK, FEAST and PRIMME.
\end{remark}

\begin{remark}
In the existing two-level convergence theory \cite{y-01,d-02},
fine and coarse time-grid propagators are assumed to be unitarily
diagonalizable (or, equivalently, normal matrices). However,
this assumption is extremely strong and it probably isn't always true,
even for \cite[Example 4.1.2]{d-02}, where the linear operator $G$
usually cannot be unitarily diagonalizable.
\end{remark}

\subsection{Specific application}

For convenience of analysis we rewrite the time-marching loop \eqref{eq-08} as,
ignoring the auxiliary vectors $\{\mathcal{V}_j\}_{j=0}^N$:
\begin{equation}\label{eq-12}
(M_s + \tau_{k+1} Q_s) \mathcal{U}_{k+1} = M_s \mathcal{U}_k + \tau_{k+1} \mathcal{F}_{k+1},~~k=0,\cdots,N-1,
\end{equation}
where the Schur complement
\begin{equation*}
Q_s=(m_{11}^zA_s+a_{11}^zM_s)-(\tilde{m}_z\otimes A_s+\tilde{a}_z\otimes M_s)
(\tilde{M}_z\otimes A_s+\tilde{A}_z\otimes M_s)^{-1}(\tilde{m}_z^T\otimes A_s+\tilde{a}_z^T\otimes M_s)
\end{equation*}
with vectors $\tilde{m}_z=(m_{12}^z,\cdots,m_{1n}^z)$ and $\tilde{a}_z=(a^z_{12},\cdots,a^z_{1n})$,
as well as submatrices $\tilde{M}_z = (m^z_{ij})_{i,j=2}^{M-1}$ and $\tilde{A}_z = (a^z_{ij})_{i,j=2}^{M-1}$.
The invertibility of $\tilde{M}_z\otimes A_s+\tilde{A}_z\otimes M_s$ is assured due to its positive definiteness,
see \cite[Proposition 2.1]{g-09}.

Regarding \eqref{eq-12}, the vector $\mathcal{G}_j = \tau_j(M_s + \tau_j Q_s)^{-1}\mathcal{F}_j$,
the fine time-grid time-stepping scheme
\begin{eqnarray}\label{eq-06}
  \Psi_j = (M_s + \tau_j Q_s)^{-1} M_s,~~j=1,\cdots,N
\end{eqnarray}
and the coarse-grid time-stepping scheme is given by 
\begin{eqnarray}\label{eq-09}
  \Phi_k = (M_s + \tilde{\tau}_k Q_s)^{-1} M_s,~~k=1,\cdots,N_c.
\end{eqnarray}
It should be emphasized that the aforementioned matrix inversion corresponds to a spatial solve,
however, not by an appropriate solver for \eqref{eq-12} but by the multigrid method with line smoother for \eqref{eq-08}.

\begin{theorem}\label{thm-02}
 Time-grid propagators \eqref{eq-06} and \eqref{eq-09} are both strongly stable.
\end{theorem}

\begin{Proof}
Viewed from the specific expressions \eqref{eq-06} and \eqref{eq-09},
it is clear that $\Psi_j$ ($j=1,\cdots,N$) and $\Phi_k$ ($k=1,\cdots,N_c$)
can be simultaneously diagonalized
by an appropriate matrix $\mathcal{X}_s$ which makes $M_s^{1/2}\mathcal{X}_s$ unitary. The reason is that
\begin{eqnarray*}
 \Psi_j = M_s^{-1/2}(I + \tau_jM_s^{-1/2} Q_s M_s^{-1/2})^{-1}M_s^{1/2}~~
 \Rightarrow ~~(I + \tau_jM_s^{-1/2} Q_s M_s^{-1/2})^{-1}=M_s^{1/2}\Psi_jM_s^{-1/2},
\end{eqnarray*}
which can be of course diagonalizable by the orthogonal matrix $M_s^{1/2}\mathcal{X}_s$
because it is symmetric (real) to ensure its normality.
It is important to notice that since matrices
$(I + \tau_jM_s^{-1/2} Q_s M_s^{-1/2})^{-1}$ and
$(I + \tilde{\tau}_kM_s^{-1/2} Q_s M_s^{-1/2})^{-1}$ are both SPD and respectively similar to $\Psi_j$ and $\Phi_k$,
the $\omega$-th eigenvalues of $\Psi_j$ and $\Phi_k$ is calculated by
\begin{eqnarray}\label{eq-13}
 \lambda^{(\omega)}_j=\frac{1}{1+\tau_j\sigma_\omega}~~\mbox{and}~~\mu^{(\omega)}_k=\frac{1}{1+\tilde{\tau}_k\sigma_\omega},~\omega=1,\cdots,n,
\end{eqnarray}
where $\sigma_\omega$ is the $\omega$-th eigenvalue of the matrix $M_s^{-1/2} Q_s M_s^{-1/2}$.
It is obvious that, for any eigenmode $\omega$,
eigenvalues $\lambda^{(\omega)}_j$ $(j=1,\cdots,N)$ and $\mu^{(\omega)}_k$ $(k=1,\cdots,N_c)$
are all real, positive and smaller than one due to the SPD property of $M_s^{-1/2} Q_s M_s^{-1/2}$
which makes $\sigma_\omega>0$ for $\omega=1,\cdots,n$. This completes the proof.
\end{Proof}

Theorem \ref{thm-02} indicates that the two key assumptions in
Stage {\bf S1} are satisfied, then, plugging \eqref{eq-13}, \eqref{eq-18} becomes
\begin{eqnarray*}
\begin{bmatrix}
 \zeta_2^{(\omega)} & -\rho_2^{(\omega)} &  &  &  &  \\
 -\rho_2^{(\omega)} & \theta_2^{(\omega)}+\zeta_3^{(\omega)} & -\rho_3^{(\omega)} &  &  &  \\
         & \ddots & \ddots & \ddots &  &  \\
 &  &  & -\rho_{N_c-2}^{(\omega)} & \theta_{N_c-2}^{(\omega)}+\zeta_{N_c-1}^{(\omega)} & -\rho_{N_c-1}^{(\omega)} \\
 &  &  &  & -\rho_{N_c-1}^{(\omega)} & \theta_{N_c-1}^{(\omega)}+\zeta_{N_c}^{(\omega)}
\end{bmatrix}~~\mbox{with}~~
\begin{matrix}
\zeta_\epsilon^{(\omega)}=\left[\pi_{\epsilon-1}^{(\omega)}\pi_\epsilon^{(\omega)}
(1+\tilde{\tau}_\epsilon\sigma_\omega)\tilde{\pi}_\epsilon^{(\omega)}\right]^2,\\
\\
\rho_\varsigma^{(\omega)}=\pi_{\varsigma-1}^{(\omega)}(\pi_\varsigma^{(\omega)})^3
(1+\tilde{\tau}_\varsigma\sigma_\omega)(\tilde{\pi}_\varsigma^{(\omega)})^2,\\
\\
\theta_\varsigma^{(\omega)}=(\pi_\varsigma^{(\omega)})^4(\tilde{\pi}_\varsigma^{(\omega)})^2,
\end{matrix}
\end{eqnarray*}
and the upper bound \eqref{eq-17} becomes
\begin{eqnarray*}
\big{\|}\mathcal{T}_\Psi^\Phi\big{\|}^2_{(\tilde{\mathcal{X}}_s\tilde{\mathcal{X}}_s^*)^{-1}}\le\max_\omega\max\left\{\delta_2^{(\omega)},
\max_{\varsigma=2,\cdots,N_c-1}\frac{\delta_{\varsigma+1}^{(\omega)}\delta_\varsigma^{(\omega)}
(1+\tilde{\tau}_\varsigma\sigma_\omega)\pi_{\varsigma-1}^{(\omega)}}{\delta_\varsigma^{(\omega)}
(1+\tilde{\tau}_\varsigma\sigma_\omega)\pi_{\varsigma-1}^{(\omega)}-
\delta_{\varsigma+1}^{(\omega)}\pi_\varsigma^{(\omega)}}\right\},
\end{eqnarray*}
where
\begin{eqnarray*}
 \pi_\varsigma^{(\omega)}=\prod_{j=(\varsigma-1)m+1}^{\varsigma m}(1+\tau_j\sigma_\omega),~~
 \tilde{\pi}_\varsigma^{(\omega)}=\frac{1}{1+\tilde{\tau}_\varsigma\sigma_\omega-\pi_\varsigma^{(\omega)}}.
\end{eqnarray*}

\section{Numerical experimentations}\label{sec4}

Numerical results are reported to illustrate the convergence rate
and the error-norms decay property of the fully discrete FE scheme \eqref{eq-04},
profile how robust should the multigrid with line smoother proposed in \cite{c-01} be
in regards to the fractional order $\alpha$ and the spatial mesh size,
and evaluate how sharp is our derived convergence upper bound listed in Theorem \ref{thm-01},
for solving problem \eqref{eq-01} posed on a unit square and an L-shaped domains.
All our numerical computations are carried out on a 64-bit Fedora 14 platform,
in double precision arithmetic on Intel Xeon (W5590) with 24.0 gigabytes random access memory, 
3.33 gigahertz and an -O2 compiler optimization option.

\subsection{{\color{red}Convergence order and robustness tests}}

\begin{example}\label{exm-01}
  Consider problem \eqref{eq-01} with {\color{red}$d=2$}, the reference solution $u=e^{-t}\sin(\pi x)\sin(\pi y)$ and $\Omega=(0,1)^2$.
\end{example}

We partition the solution domain $\Omega$ into a uniform triangular mesh with size $h_{\Omega}$,
and construct a graded mesh in the extended direction using the formula \eqref{eq-07} with the length $\mathcal{Z}=1$.
The finite element space and respective matrices were introduced in Section \ref{sec2}.
The chosen time-step sizes $\tau = 10^{-4}$ and $\tau = 5\times10^{-5}$.
The tolerance for stopping of the multigrid V(1,1)-cycle algorithm
(using the zero initial guess and one presmoothing and postsmoothing sweep
performed by the vertical line smoother in red-black ordering) in solving linear systems is
that the discrete $l^2$-norm of the residual is reduced by a factor of $10^8$ relative to that of the initial residual.
Table \ref{tlb-01} shows the number of degrees of freedom (DoF),
the error-norms $$\|e\|^\tau_N:=\|u^\mathcal{J}_N - u(x,y,N\tau)\|_{H^1(\Omega)},~~
\|e\|^\tau_{\max}:=\max_{j=1,\cdots,N}\|u^\mathcal{J}_j-u(x,y,j\tau)\|_{H^1(\Omega)},$$
the respective rate computed by the error of the coarser grid divided by that of the current grid,
and the average number of iterations for multigrid
required to achieve convergence for $\alpha=0.4$, $1.0$ and $1.4$. In these cases, the results are all in agreement with theory,
and the multigrid V(1,1)-cycle is almost uniform with respect to the fractional order $\alpha$ and the number of DoF
(or, equivalently, the spatial mesh size $h_{\Omega}$).
In particular, we can observe that the computational
error-norm decay rate in $\|\tilde{\mathcal{U}}-U^\mathcal{J}_N\|_{\mathring{H}^1_L(z^\beta,\mathcal{C})}$ is about DoF$^{-1/3}$.

\begin{table}[htbp]
\small
\begin{center}
\caption{Error results, convergence rates and average number of iterations for three different $\alpha$ and two different $\tau$.}\label{tlb-01}\vskip 0.1cm
\begin{tabular}{ccc|ccccc|ccccc}\hline
$\alpha$ & $h_{\Omega}$ & DoF & $\|e\|^{1\textrm{e}-4}_{100}$ & rate & $\|e\|^{1\textrm{e}-4}_{\max}$ & rate & aIter
& $\|e\|^{5\textrm{e}-5}_{200}$ & rate & $\|e\|^{5\textrm{e}-5}_{\max}$ & rate & aIter \\ \hline
\multirow{5}{*}{0.4}
&1/4      &   225   & 8.37e-1 &      & 8.42e-1 &      & 11 & 8.39e-1 &      & 8.42e-1 &      & 11 \\ 
&1/8      &  1620   & 4.28e-1 & 1.96 & 4.43e-1 & 1.90 & 12 & 4.29e-1 & 1.96 & 4.30e-1 & 1.96 & 12 \\ 
&1/16     &  12716  & 2.16e-1 & 1.98 & 2.17e-1 & 2.04 & 12 & 2.16e-1 & 1.99 & 2.17e-1 & 1.98 & 11 \\ 
&1/32     &  104544 & 1.09e-1 & 1.98 & 1.10e-1 & 1.97 & 12 & 1.09e-1 & 1.98 & 1.10e-1 & 1.97 & 12 \\ 
&1/64     &  887250 & 5.58e-2 & 1.95 & 5.87e-2 & 1.87 & 12 & 5.58e-2 & 1.95 & 5.87e-2 & 1.87 & 12 \\ \hline
\multirow{5}{*}{1.0}
&1/4      &   225   & 8.39e-1 &      & 8.42e-1 &      & 11 & 8.37e-1 &      & 8.42e-1 &      & 11 \\ 
&1/8      &  1620   & 4.30e-1 & 1.95 & 4.30e-1 & 1.96 & 12 & 4.28e-1 & 1.96 & 4.30e-1 & 1.96 & 12 \\ 
&1/16     &  12716  & 2.16e-1 & 1.99 & 2.17e-1 & 1.98 & 13 & 2.15e-1 & 1.99 & 2.17e-1 & 1.98 & 12 \\ 
&1/32     &  104544 & 1.09e-1 & 1.98 & 1.10e-1 & 1.97 & 13 & 1.08e-1 & 1.99 & 1.10e-1 & 1.97 & 13 \\ 
&1/64     &  887250 & 5.40e-2 & 2.02 & 5.86e-2 & 1.88 & 13 & 5.40e-2 & 2.00 & 5.86e-2 & 1.88 & 13 \\ \hline
\multirow{5}{*}{1.4}
&1/4      &   225   & 8.35e-1 &      & 8.42e-1 &      & 11 & 8.36e-1 &      & 8.42e-1 &      & 11 \\ 
&1/8      &  1620   & 4.28e-1 & 1.95 & 4.30e-1 & 1.96 & 12 & 4.28e-1 & 1.95 & 4.30e-1 & 1.96 & 12 \\ 
&1/16     &  12716  & 2.16e-1 & 1.98 & 2.17e-1 & 1.98 & 13 & 2.16e-1 & 1.98 & 2.17e-1 & 1.98 & 13 \\ 
&1/32     &  104544 & 1.09e-1 & 1.98 & 1.10e-1 & 1.97 & 13 & 1.09e-1 & 1.98 & 1.10e-1 & 1.97 & 13 \\ 
&1/64     &  887250 & 5.60e-2 & 1.95 & 5.85e-2 & 1.88 & 14 & 5.60e-2 & 1.95 & 5.86e-2 & 1.88 & 13 \\ \hline
\end{tabular}
\end{center}
\end{table}

{\color{red}Next we study the dependencies of the multigrid V(1,1)-cycle algorithm over the fractional order
$\alpha$ and time-step size $\tau$. The results refer to five different $\alpha$ and five different $\tau$
as shown in Table \ref{tlb-02}, where $h_{\Omega}=1/16$ and $T=1/4$.
It can be noticed that the multigrid V(1,1)-cycle algorithm converges robustly with respect to $\alpha$ and $\tau$.}

\begin{table}[htbp]
\small
\begin{center}
\caption{{\color{red}The average number of multigrid iterations for five different
$\alpha$ and five different $\tau$ with $h_{\Omega}=1/16$ and $T=1/4$.}}\label{tlb-02}\vskip 0.1cm
\begin{tabular}{c|cccccccccccc}\hline
~ & $\tau=1/64$ & $\tau=1/128$ & $\tau=1/256$ & $\tau=1/512$ & $\tau=1/1024$ \\ \hline
$\alpha=0.6$ & 12 & 12 & 12 & 11 & 11 \\ 
$\alpha=0.8$ & 13 & 12 & 12 & 12 & 11 \\ 
$\alpha=1.0$ & 13 & 13 & 12 & 12 & 12 \\ 
$\alpha=1.2$ & 13 & 13 & 13 & 12 & 12 \\ 
$\alpha=1.6$ & 14 & 13 & 13 & 13 & 12 \\ \hline
\end{tabular}
\end{center}
\end{table}

\begin{example}\label{exm-02}
  Consider problem \eqref{eq-01} with {\color{red}$d=2$}, the reference solution $u=e^{-t}\sin(\pi x)\sin(\pi y)$ and $\Omega=(-1,1)^2 \backslash(0,1)^2$.
\end{example}

\begin{figure}[htbp]
\centering
\subfigure[$\alpha = 0.4$: error-norm decay]{
\label{fig-03a}
\includegraphics[scale=0.333]{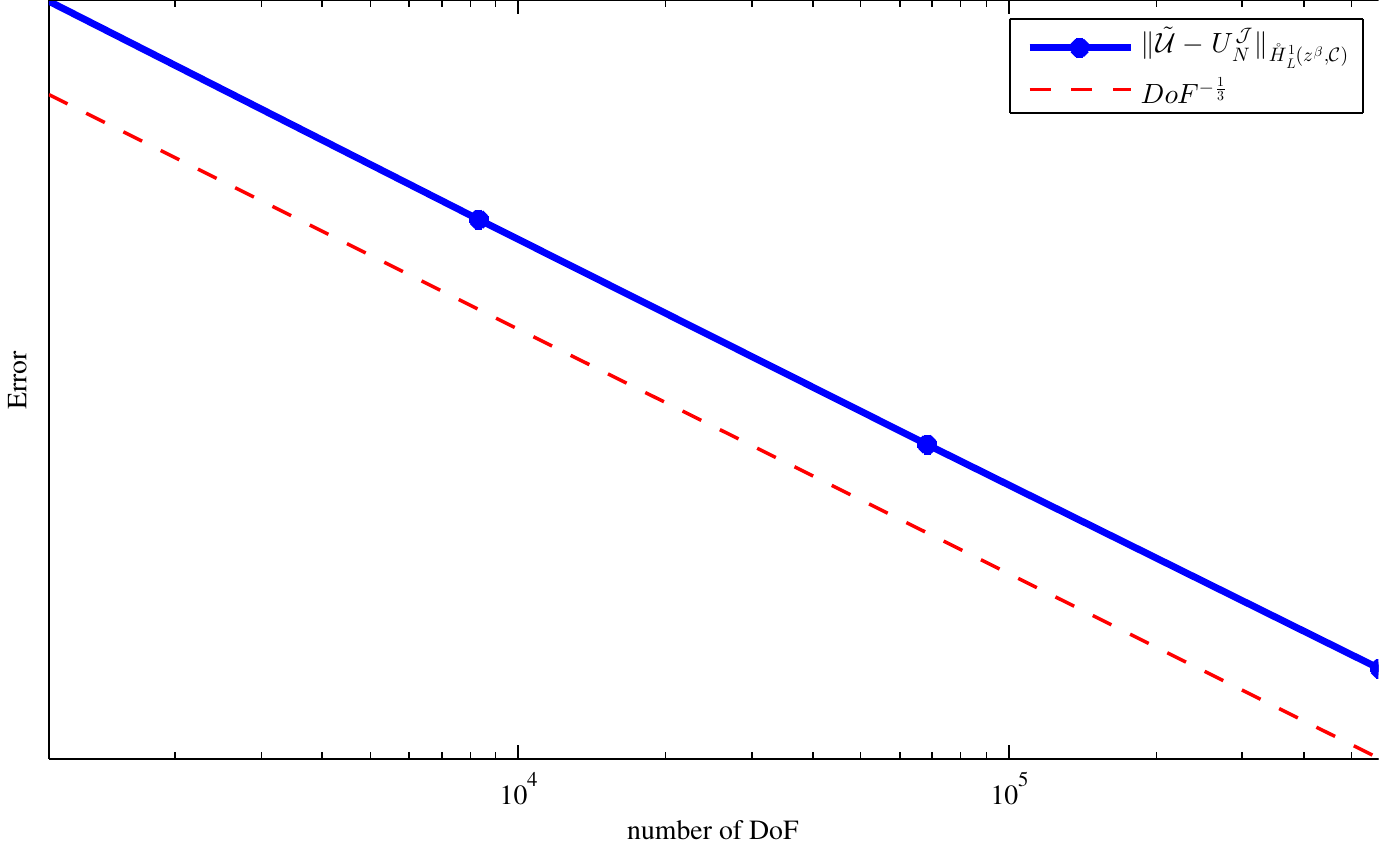}}
\hfill
\subfigure[$\alpha = 0.4$: numerical solution]{
\label{fig-03b}
\includegraphics[scale=0.333]{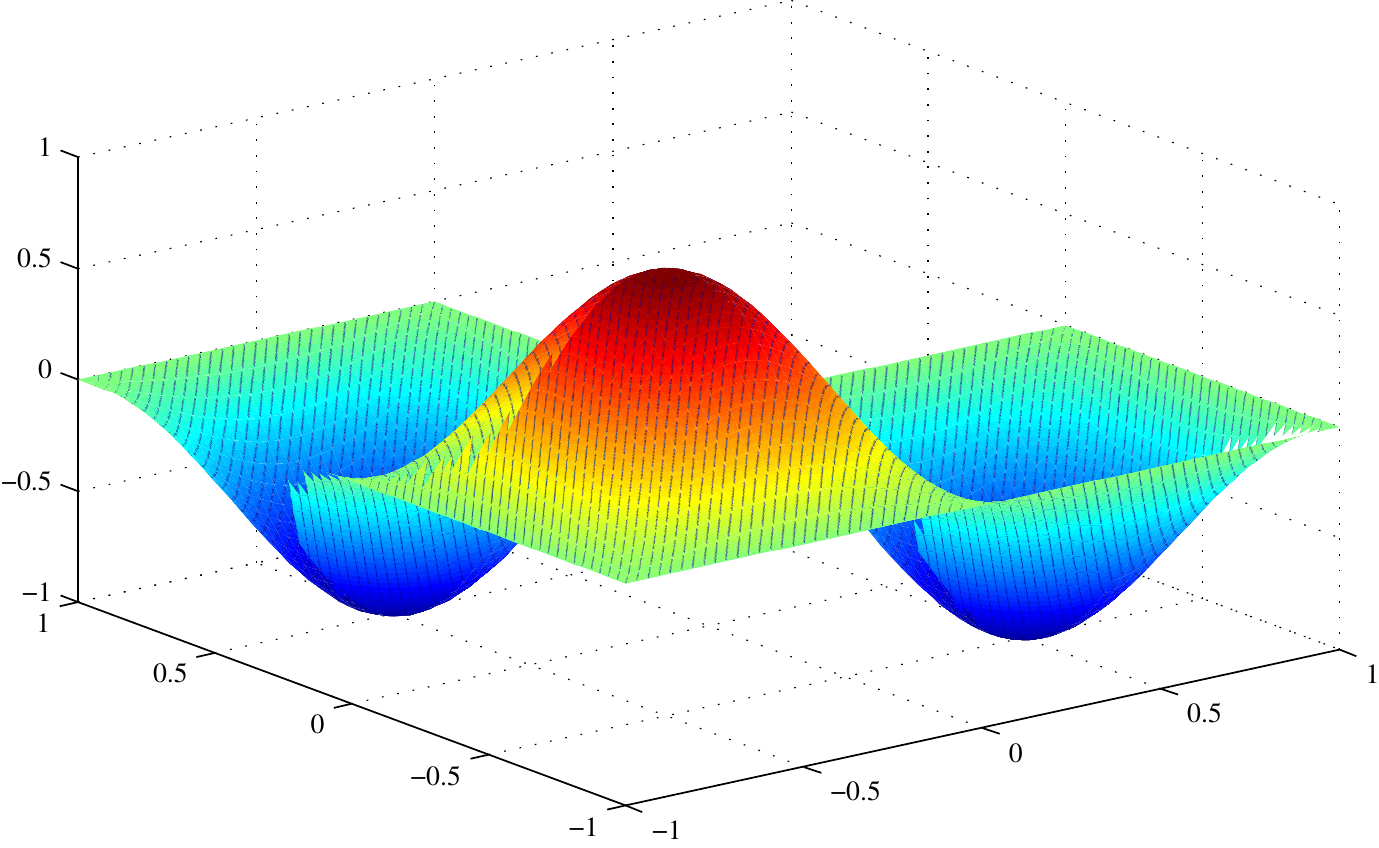}}
\hfill
\subfigure[$\alpha = 0.4$: absolute errors]{
\label{fig-03c}
\includegraphics[scale=0.333]{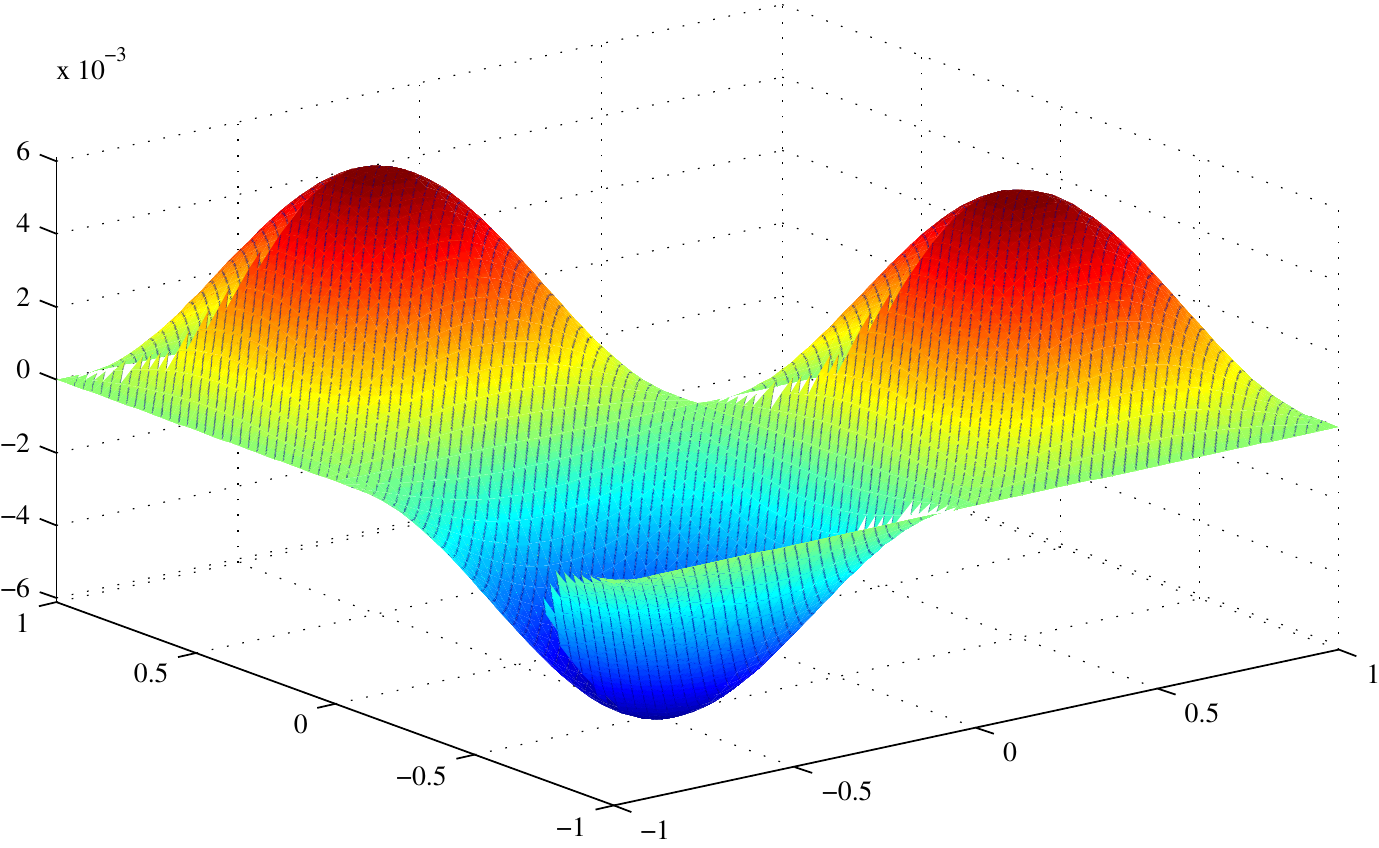}}
\centering
\subfigure[$\alpha = 1.0$: error-norm decay]{
\label{fig-03d}
\includegraphics[scale=0.333]{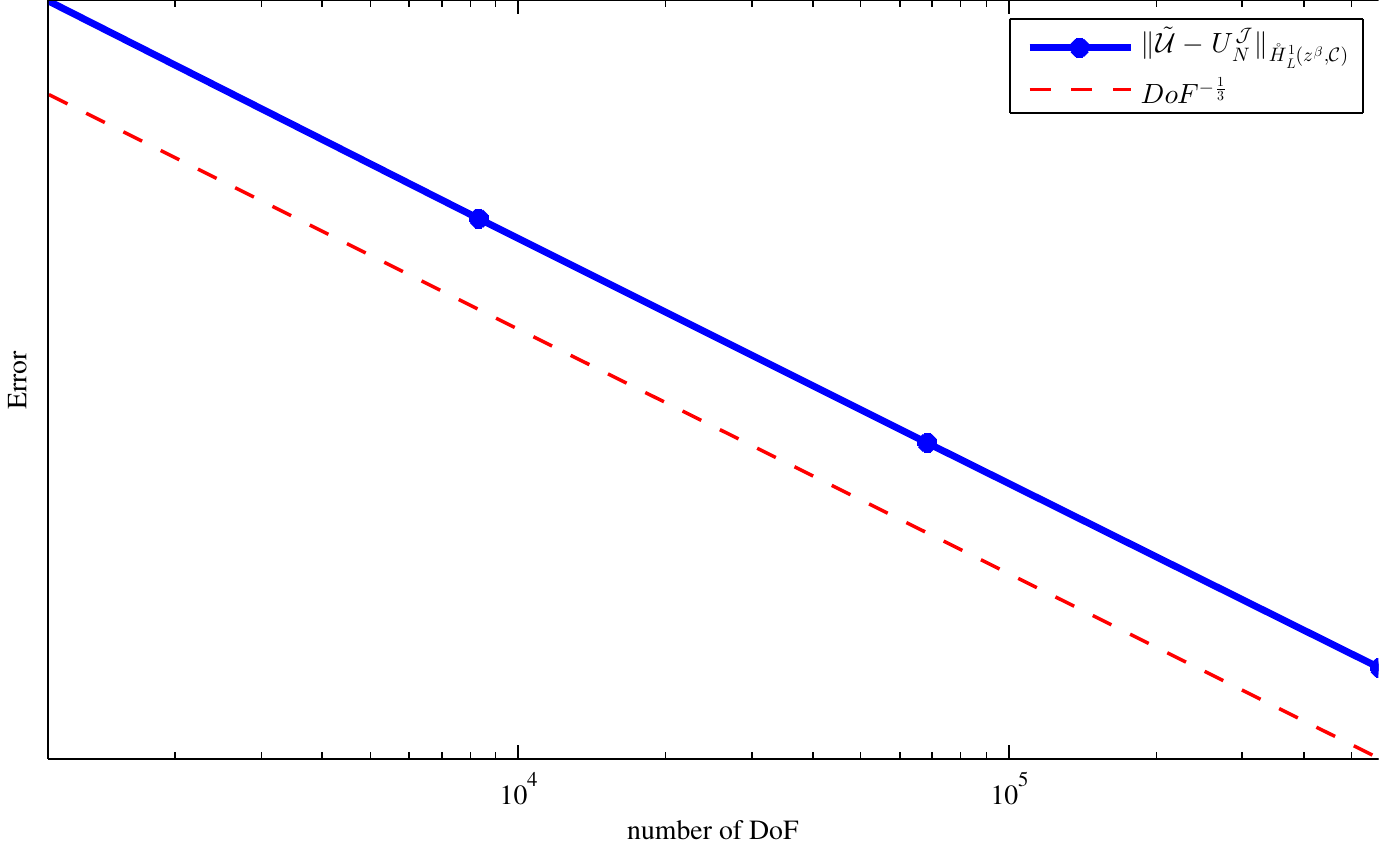}}
\hfill
\subfigure[$\alpha = 1.0$: numerical solution]{
\label{fig-03e}
\includegraphics[scale=0.333]{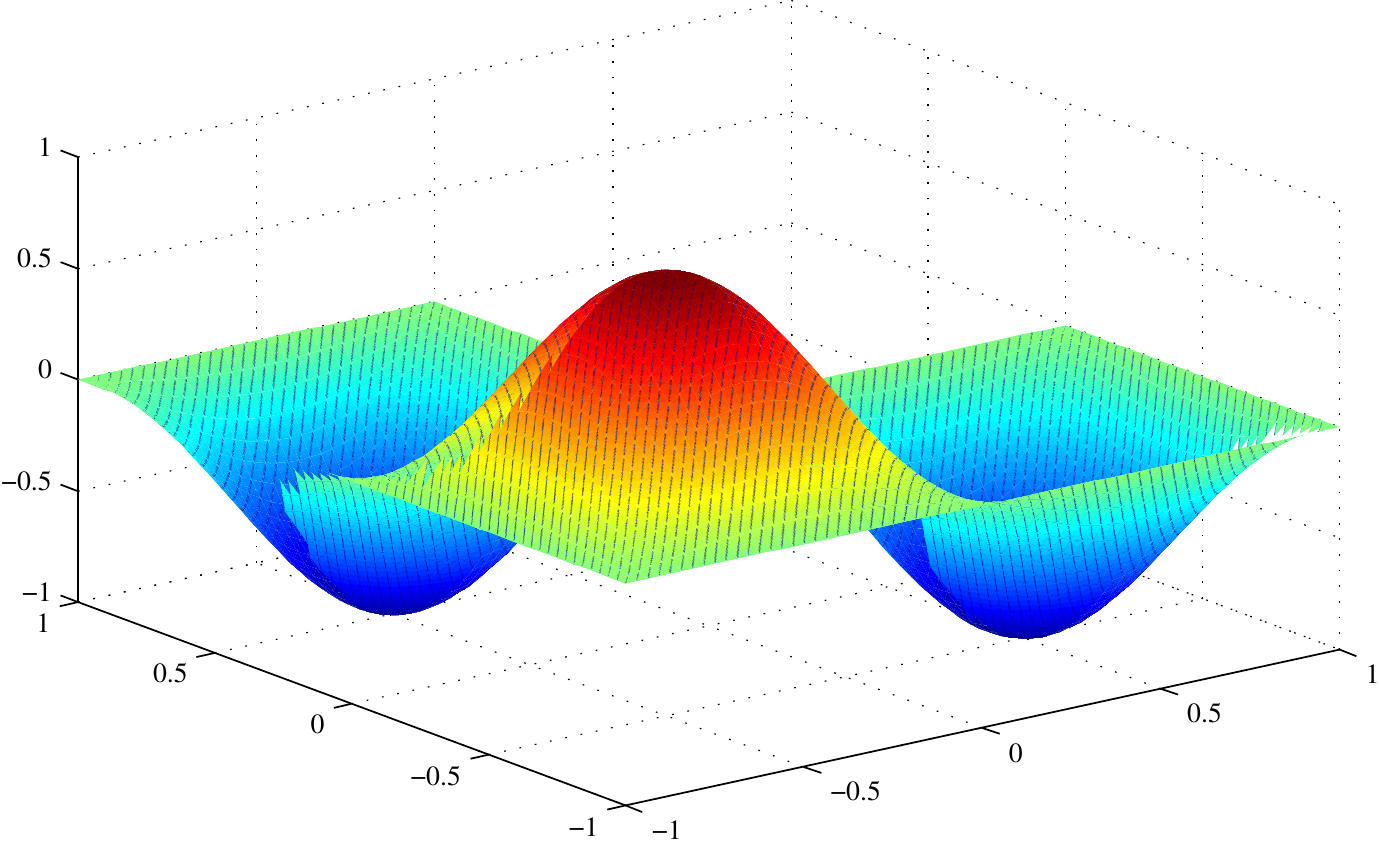}}
\hfill
\subfigure[$\alpha = 1.0$: absolute errors]{
\label{fig-03f}
\includegraphics[scale=0.333]{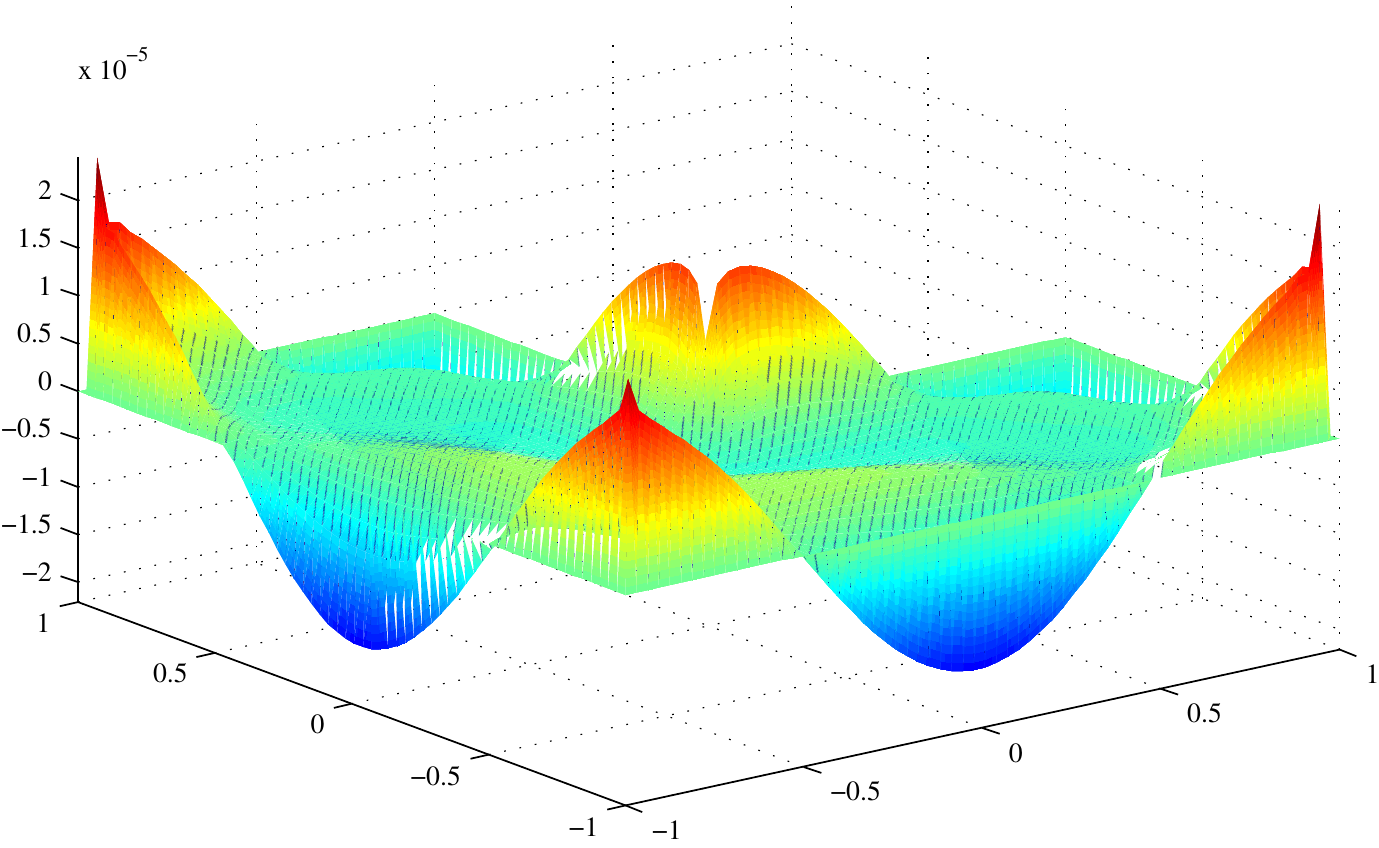}}
\centering
\subfigure[$\alpha = 1.4$: error-norm decay]{
\label{fig-03g}
\includegraphics[scale=0.333]{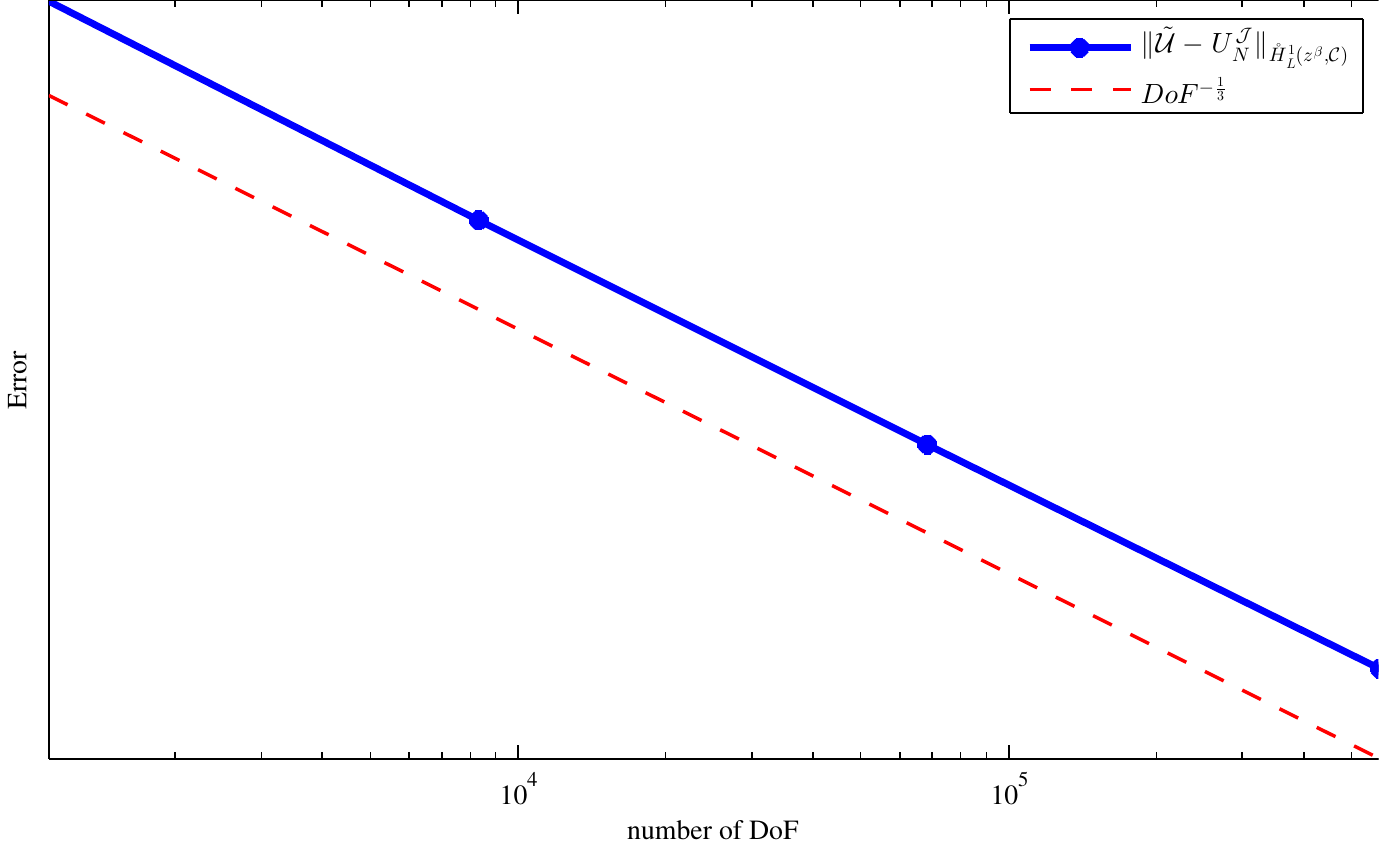}}
\hfill
\subfigure[$\alpha = 1.4$: numerical solution]{
\label{fig-03h}
\includegraphics[scale=0.333]{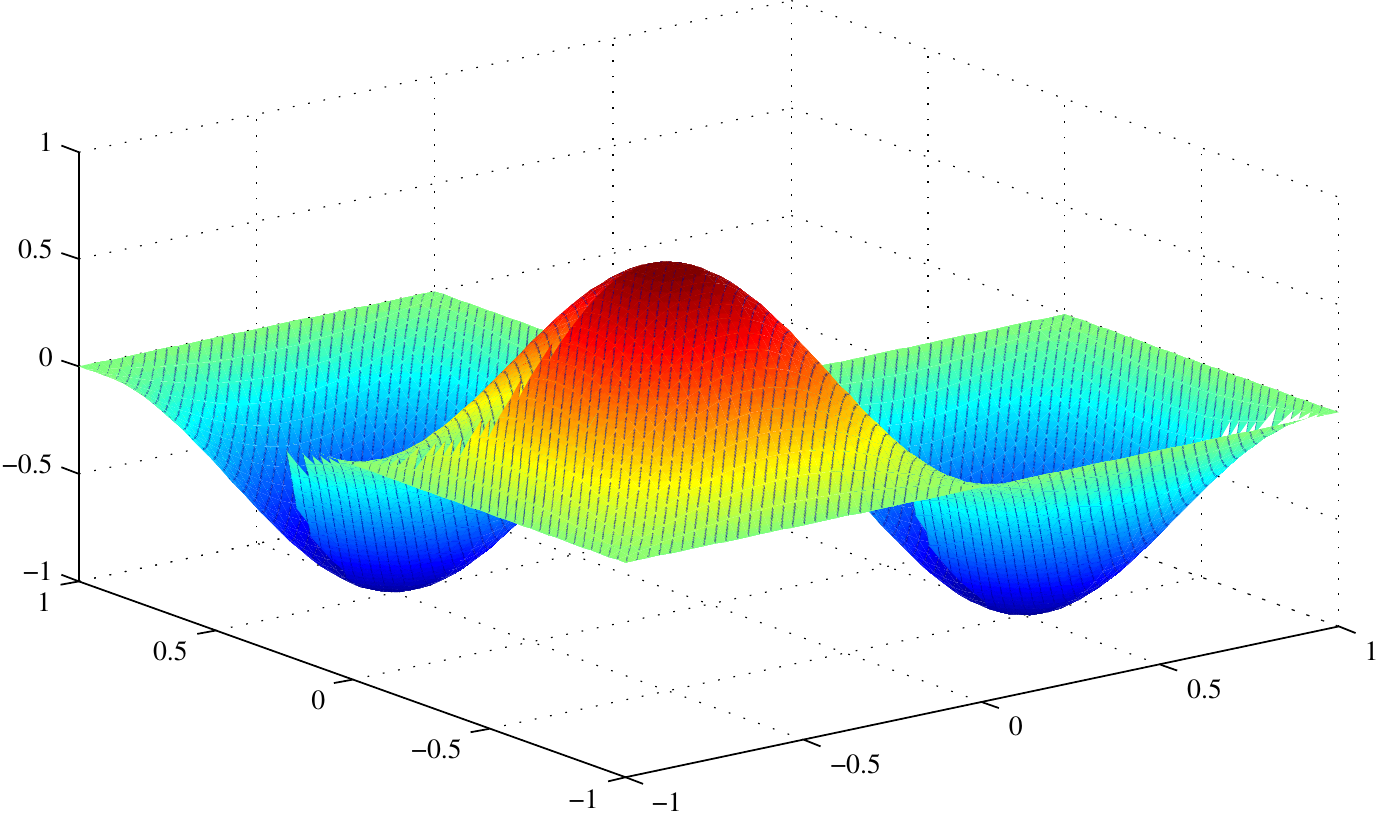}}
\hfill
\subfigure[$\alpha = 1.4$: absolute errors]{
\label{fig-03i}
\includegraphics[scale=0.333]{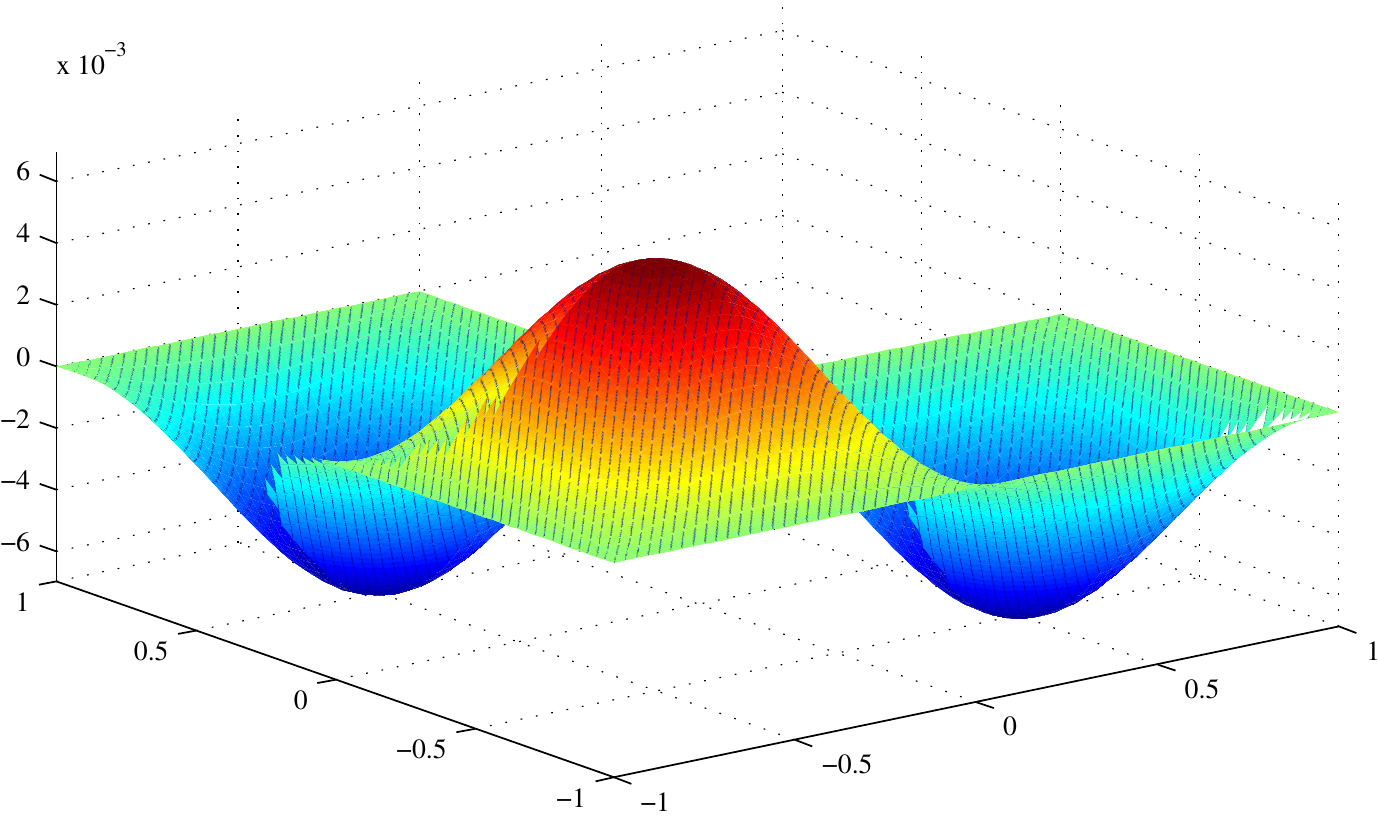}}
\caption{Computational error-norm decays, numerical solutions and absolute errors at
$T=0.01$ of Example \ref{exm-02}, 
$h_{\Omega}=1/32$, $\tau=5\times10^{-5}$.}\label{fig-03}
\end{figure}

Fig. \ref{fig-03} illustrates the computational error-norm decays, numerical solutions
and absolute errors of the reference solution and numerical solution with $\alpha = 0.4$, $1.0$ and $1.4$.
In three cases, we again obtain that the computational error-norm decay rate $\approx$ DoF$^{-1/3}$
and the fully discrete FE approximation is in much better agreement with the reference solution.

\subsection{Sharpness test in the convergence upper bound}

In this subsection, the MGRIT V(1,0)-cycle algorithm using FCF-relaxation depicted
in subsection \ref{sec3-1} together with the MGRIT V(1,0)-cycle algorithm with F-relaxation
(an equivalent version of parareal) are exploited for problem \eqref{eq-01},
where we the maximum number of time levels $L=2$.
We principally investigate the sharpness of \cite[Theorem 35]{s-03} and Theorem \ref{thm-01} in the present paper,
each of them is appended with $\sqrt{m}$ from Stage {\bf S2} in subsection \ref{sec3-2},
by a comparison with the largest observed convergence factor in terms of $\max_i\{\|r_{i+1}\|/\|r_i\|\}$,
where $r_i$ is the residual vector obtained at the $i$-th iteration.
To this end, we consider the spatial domains $\Omega=(0,1)^2$ and $\Omega=(-1,1)^2 \backslash(0,1)^2$,
and the final time $T=1$. The spatial discretization takes advantage of the uniform triangulation
with the grid spacing $h_{\Omega}=1/8$. The uniform and graded temporal partitions are involved.
The graded temporal mesh is defined by points $t_k=(k/N)^\varpi$ \cite{b-02} for the parameter $\varpi=5/2$.
The considered number of time periods $N\in\{256,1024,4096\}$ and the coarsening factor $m\in\{2,4,8,16\}$.
The space-time approximation is achieved until the discrete $l^2$ space-time residual norm
is smaller than the absolute halting tolerance $10^{-8}$, along with 100 as the maximum number of MGRIT iterations.

All test cases are implemented in C, using the open-source packages XBraid and LAPACK.
Furthermore, random initial guesses is selected for the entire temporal grid hierarchy
beyond the initial condition in \eqref{eq-01} is used at $t=0$ on the finest grid.

\begin{figure}[htbp]
\centering
\subfigure[unit square domain, $N=256$]{
\label{fig-04a}
\includegraphics[scale=0.333]{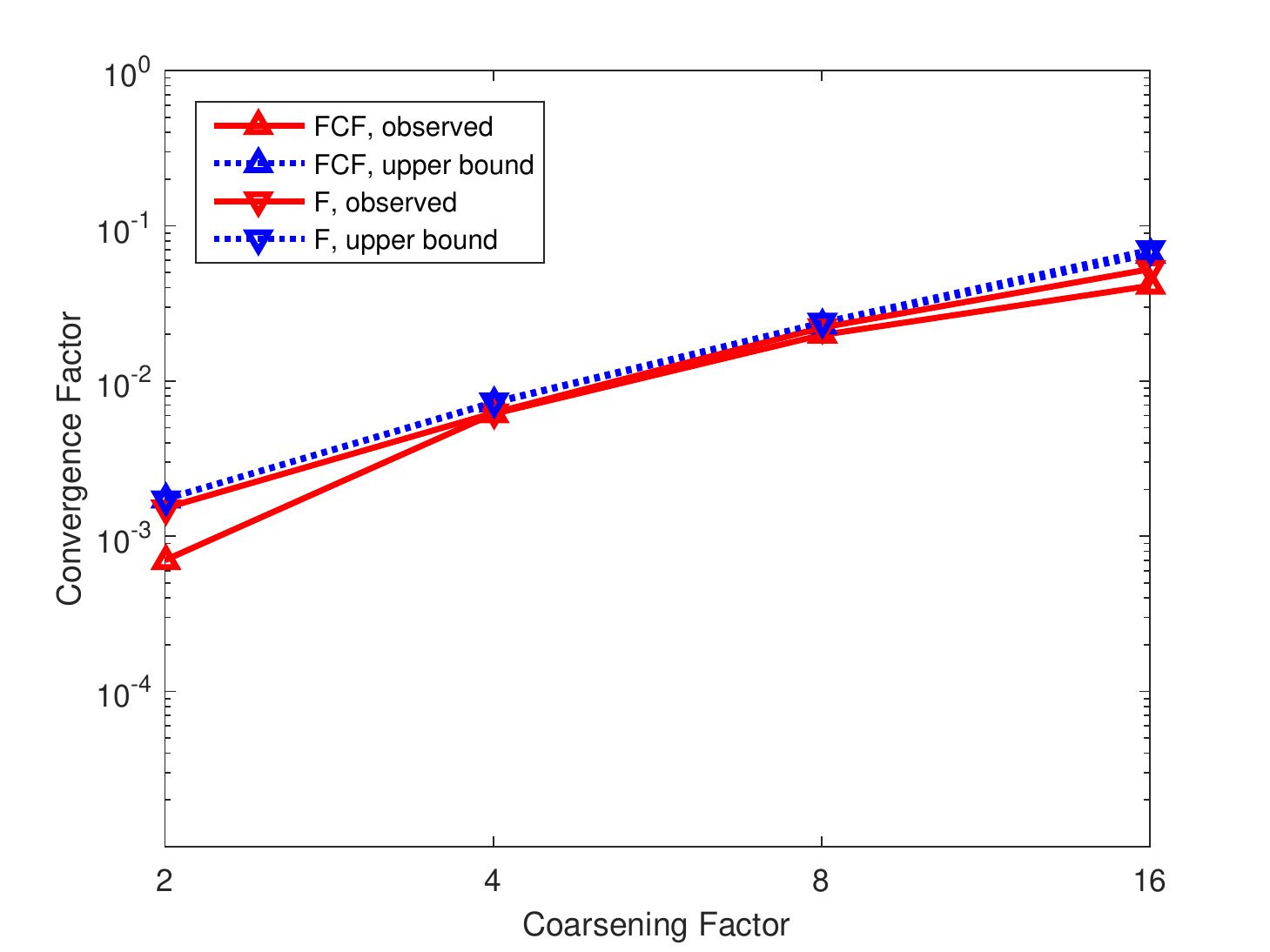}}
\hfill
\subfigure[unit square domain, $N=1024$]{
\label{fig-04b}
\includegraphics[scale=0.333]{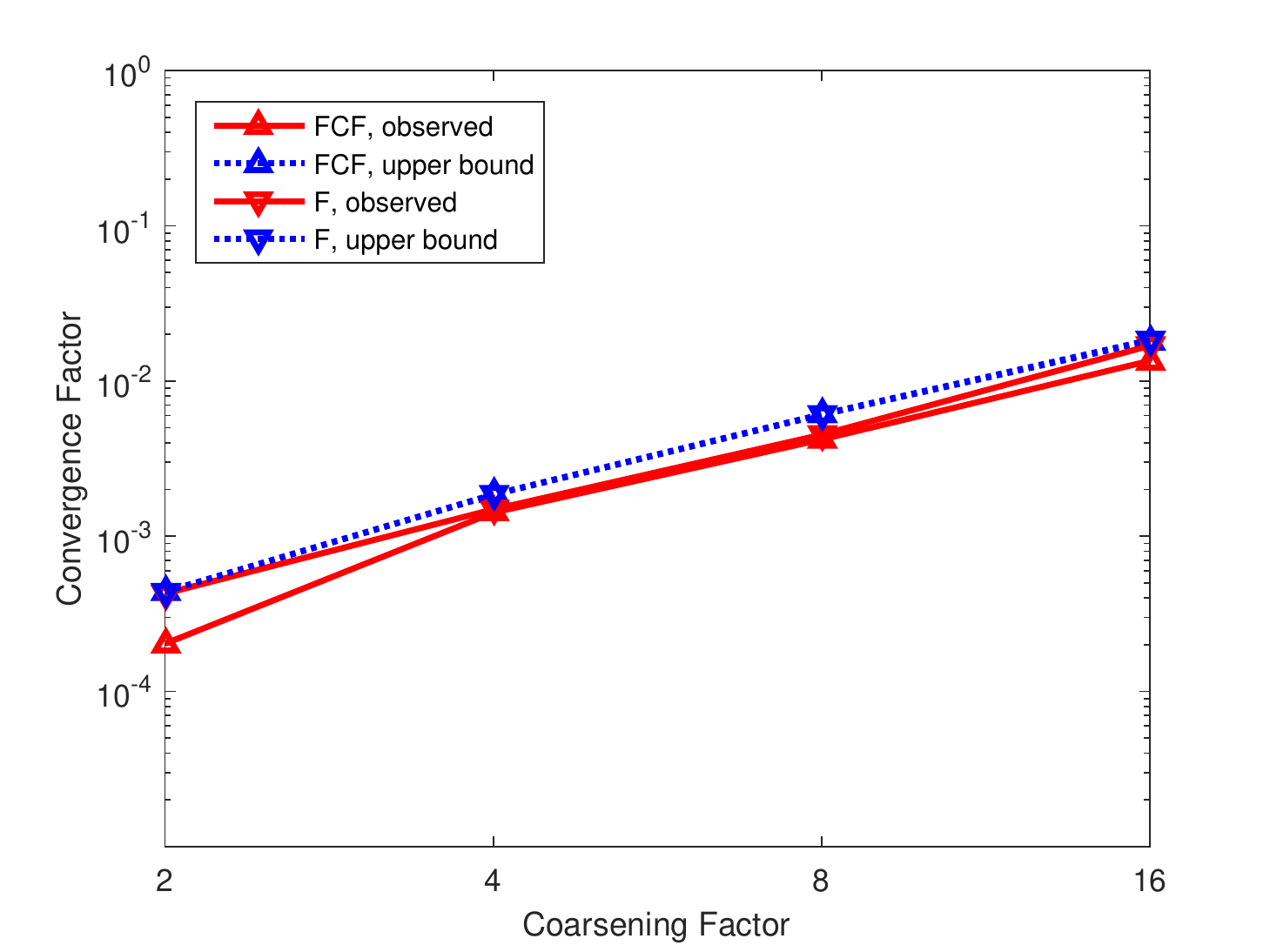}}
\hfill
\subfigure[unit square domain, $N=4096$]{
\label{fig-04c}
\includegraphics[scale=0.333]{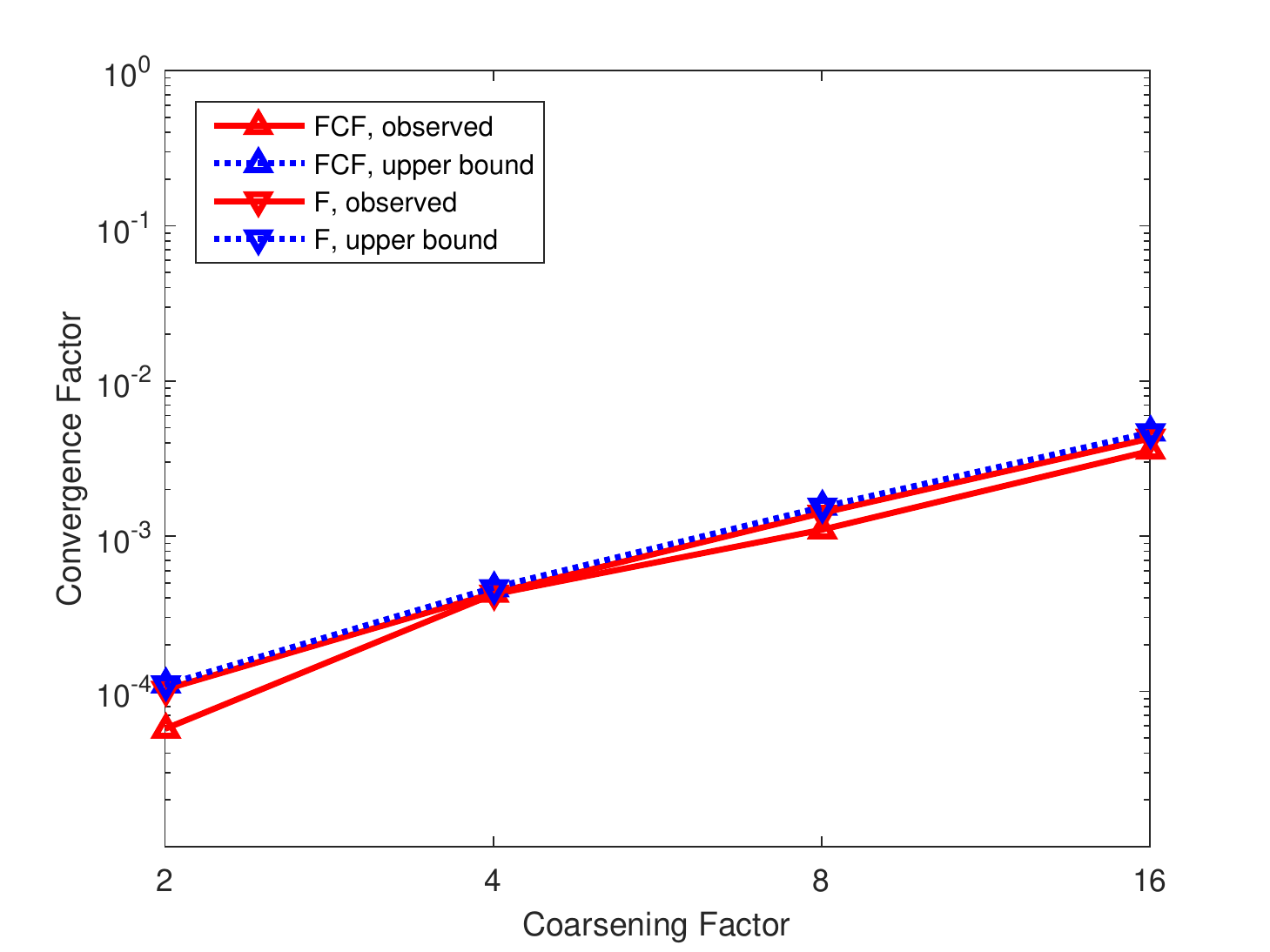}}
\centering
\subfigure[L-shaped domain, $N=256$]{
\label{fig-04d}
\includegraphics[scale=0.333]{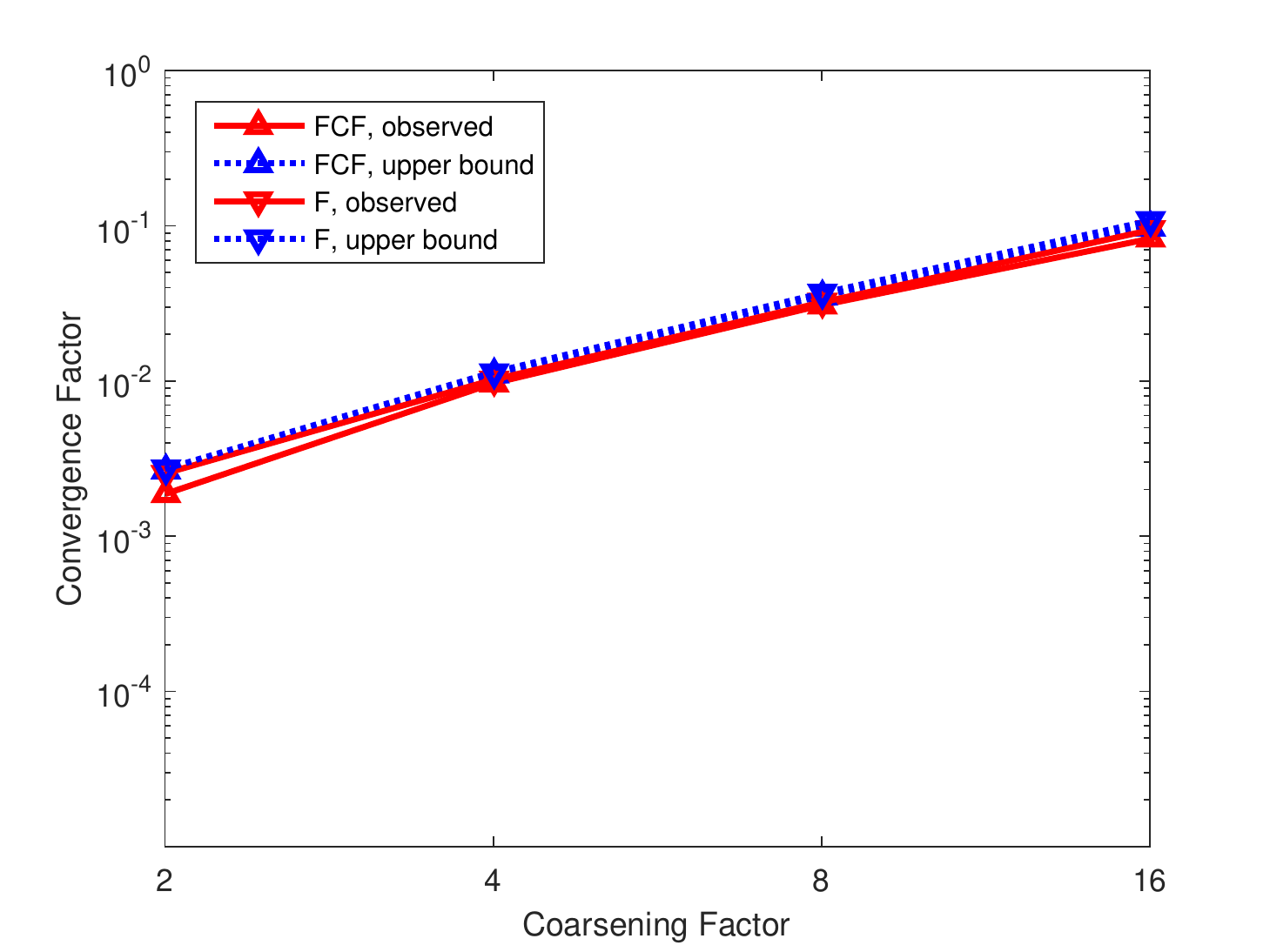}}
\hfill
\subfigure[L-shaped domain, $N=1024$]{
\label{fig-04e}
\includegraphics[scale=0.333]{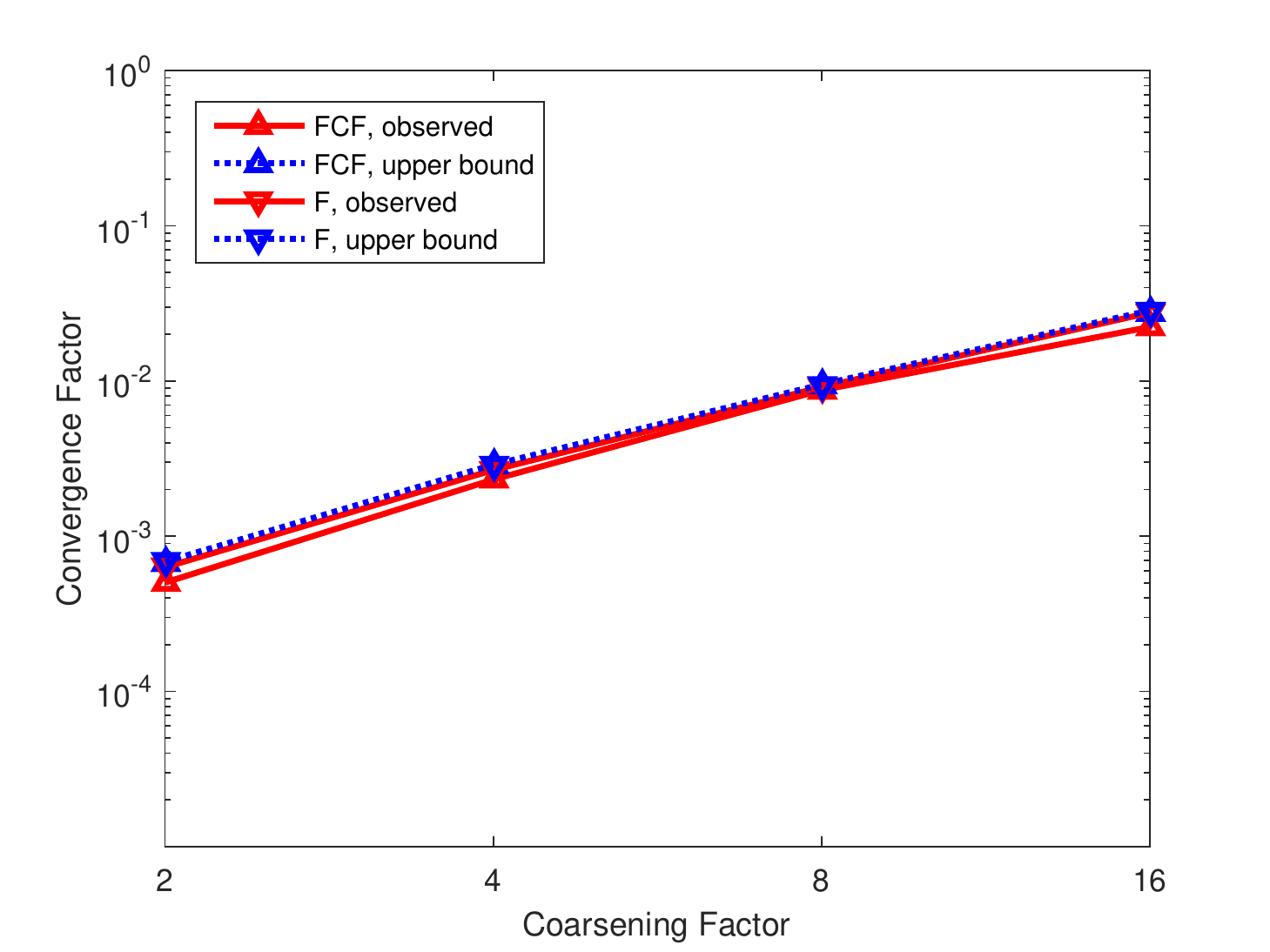}}
\hfill
\subfigure[L-shaped domain, $N=4096$]{
\label{fig-04f}
\includegraphics[scale=0.333]{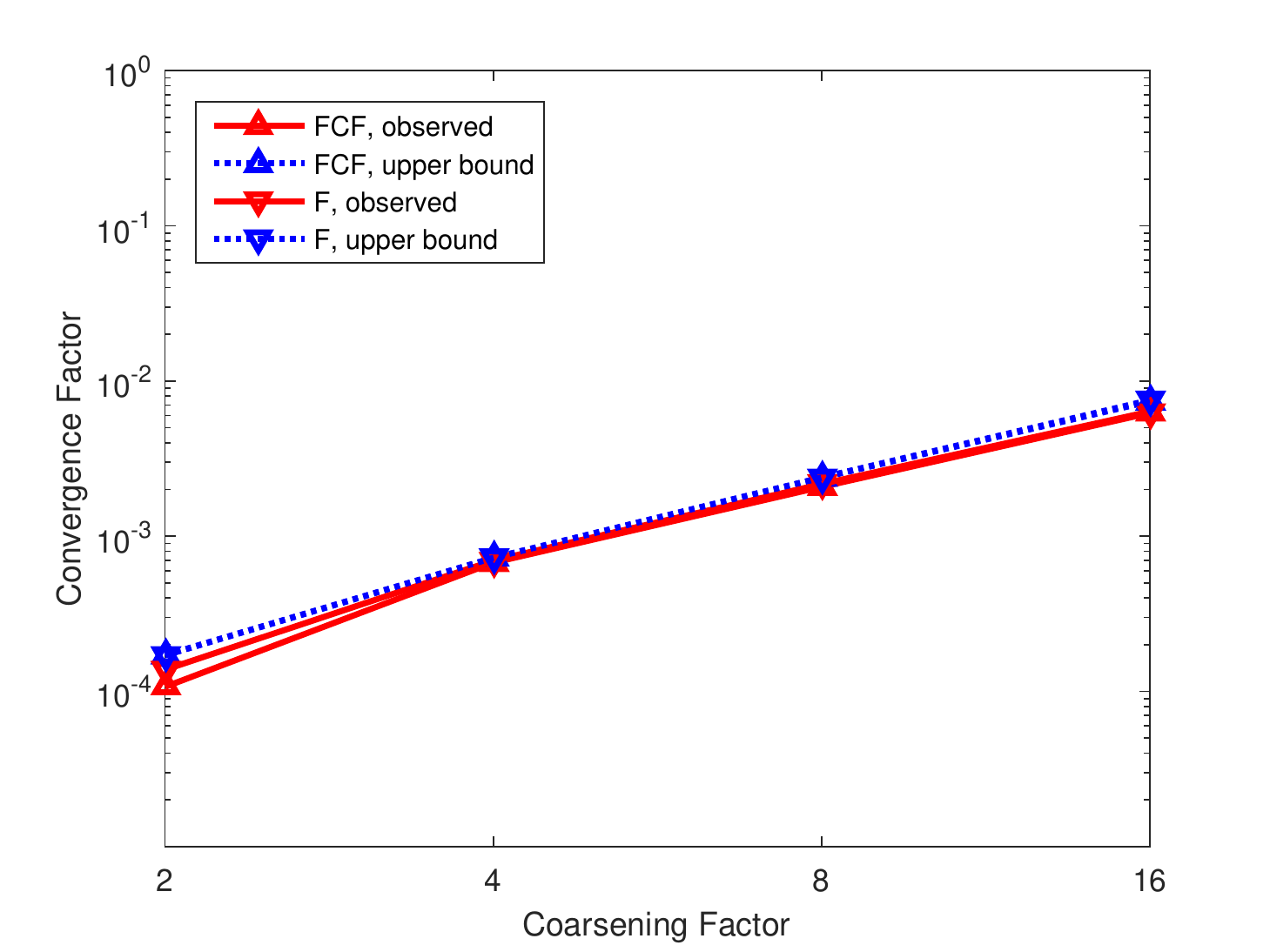}}
\caption{Uniform temporal partition: comparison on the convergence behavior of the
two-level MGRIT V(1,0)-cycle with F- and FCF-relaxation.}\label{fig-04}
\end{figure}

Results for uniform temporal partition are shown in Fig. \ref{fig-04}, demonstrating that,
in this case, all observed convergence factors are smaller than and quite close to their theoretical upper bounds.
It is indicated that single-iteration upper bounds proposed in \cite[Theorem 30]{s-03}
(two particular cases of \cite[Theorem 35]{s-03} and Theorem \ref{thm-01}) are tight.
Another observation is that results for F-relaxation are slightly larger than those based on FCF-relaxation,
even though the former corresponds to an optimal two-level algorithm {\color{red}\cite{f-02,g-02}}.
Furthermore, it must be emphasized that they decrease as the number of time periods $N$ increases
while {\color{red}they} tend to become larger as the coarsening factor $m$ increases.

\begin{figure}[htbp]
\centering
\subfigure[unit square domain, $N=256$]{
\label{fig-05a}
\includegraphics[scale=0.333]{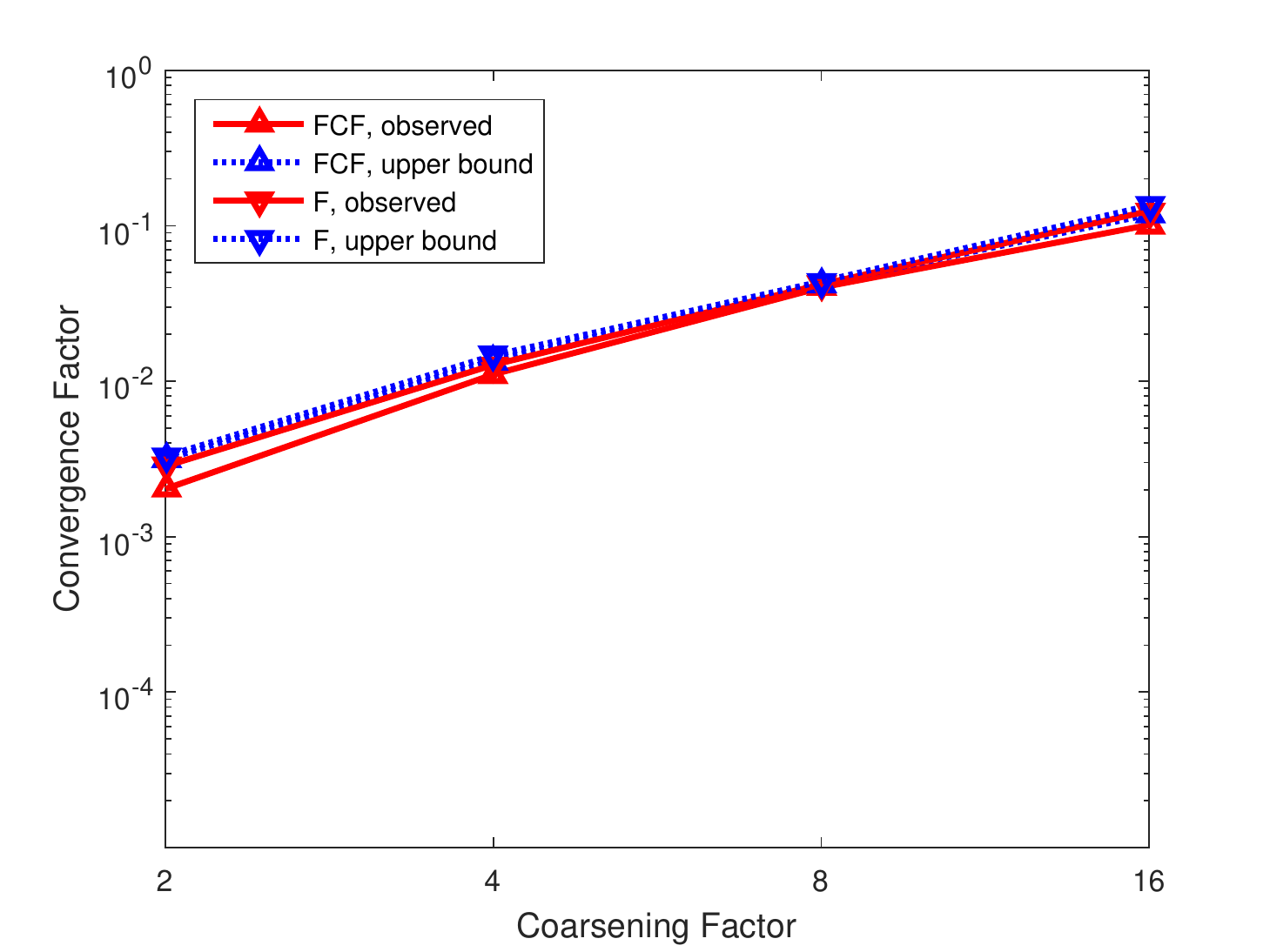}}
\hfill
\subfigure[unit square domain, $N=1024$]{
\label{fig-05b}
\includegraphics[scale=0.333]{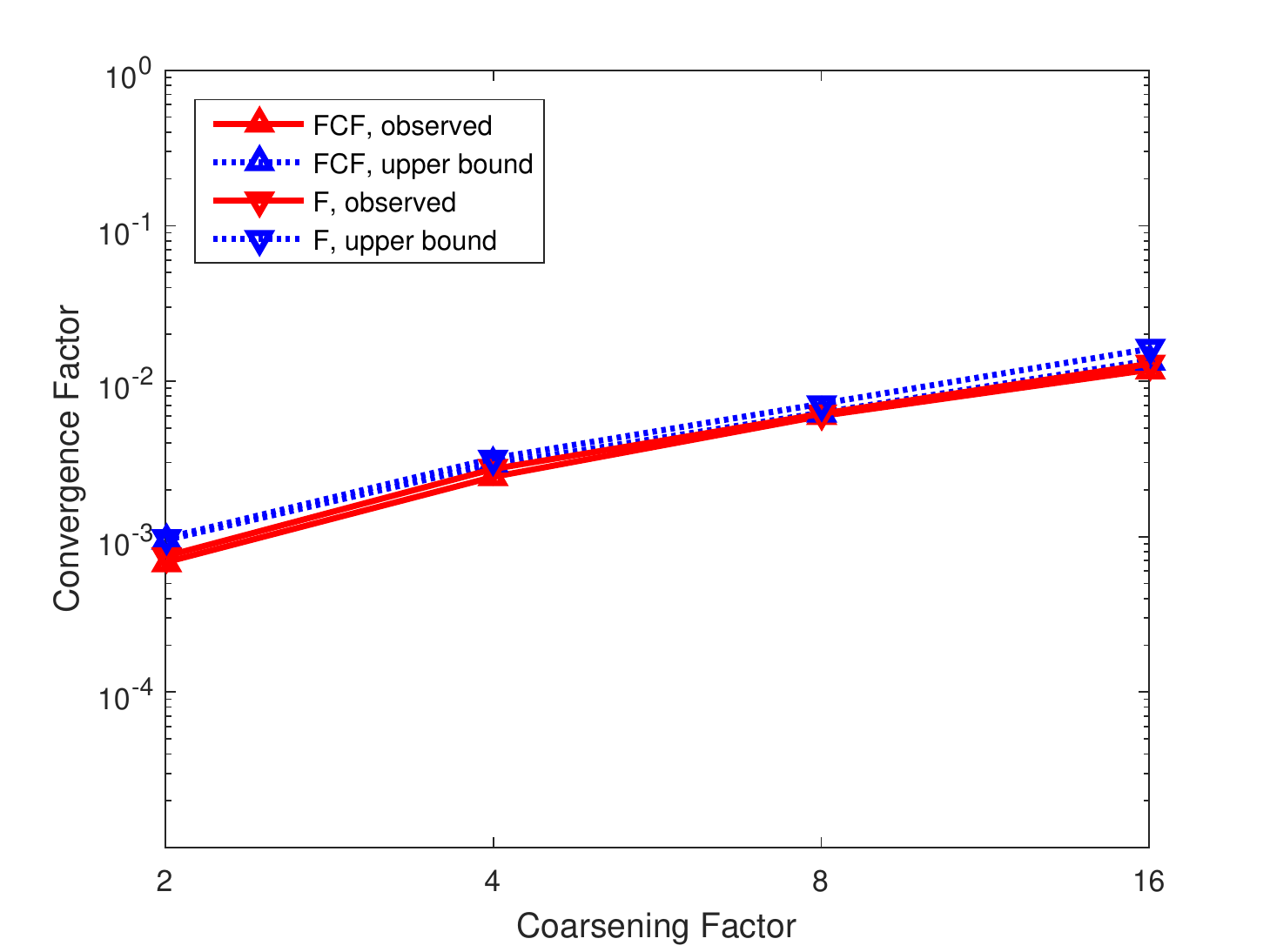}}
\hfill
\subfigure[unit square domain, $N=256$]{
\label{fig-05c}
\includegraphics[scale=0.333]{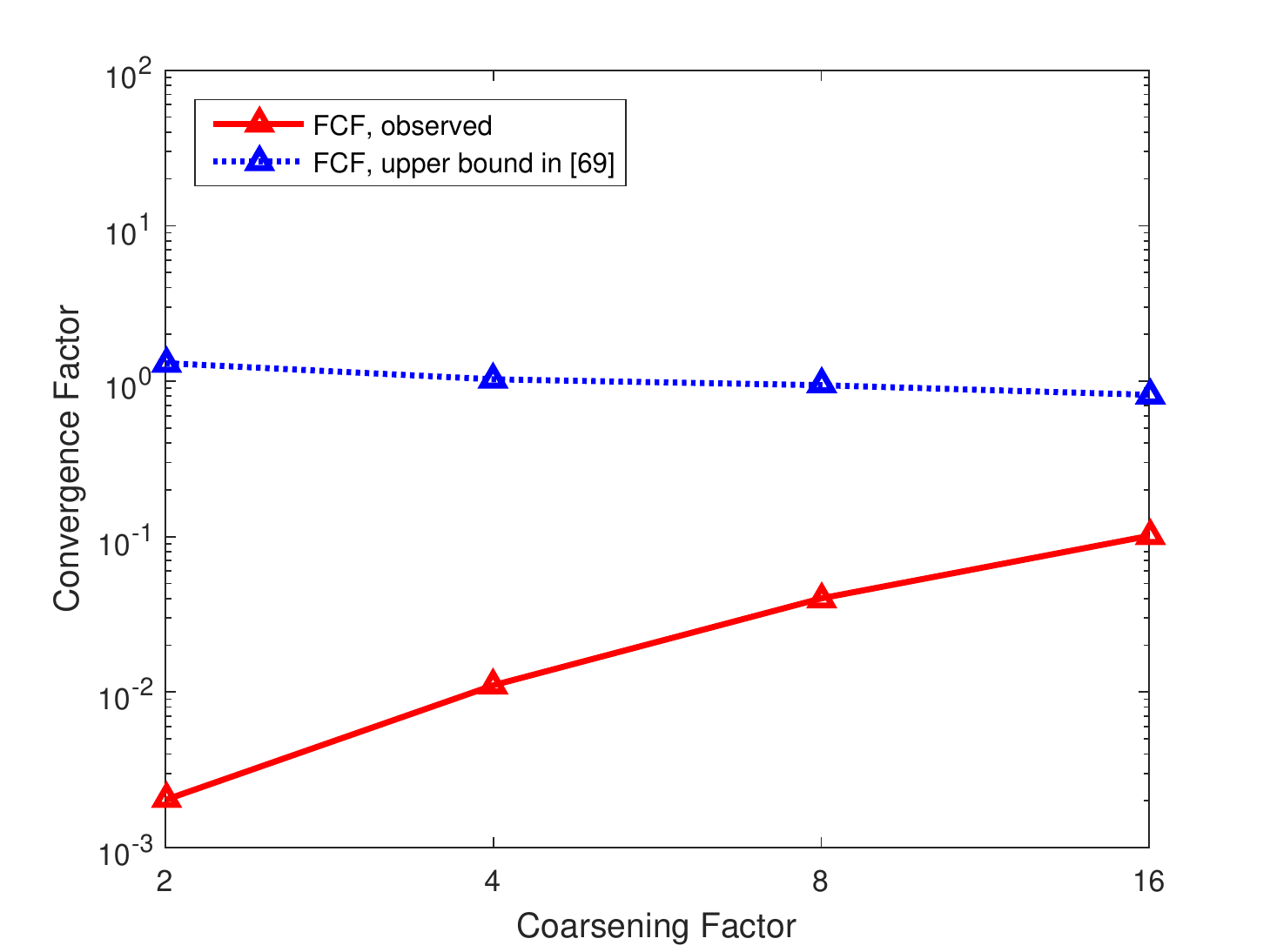}}
\centering
\subfigure[L-shaped domain, $N=256$]{
\label{fig-05d}
\includegraphics[scale=0.333]{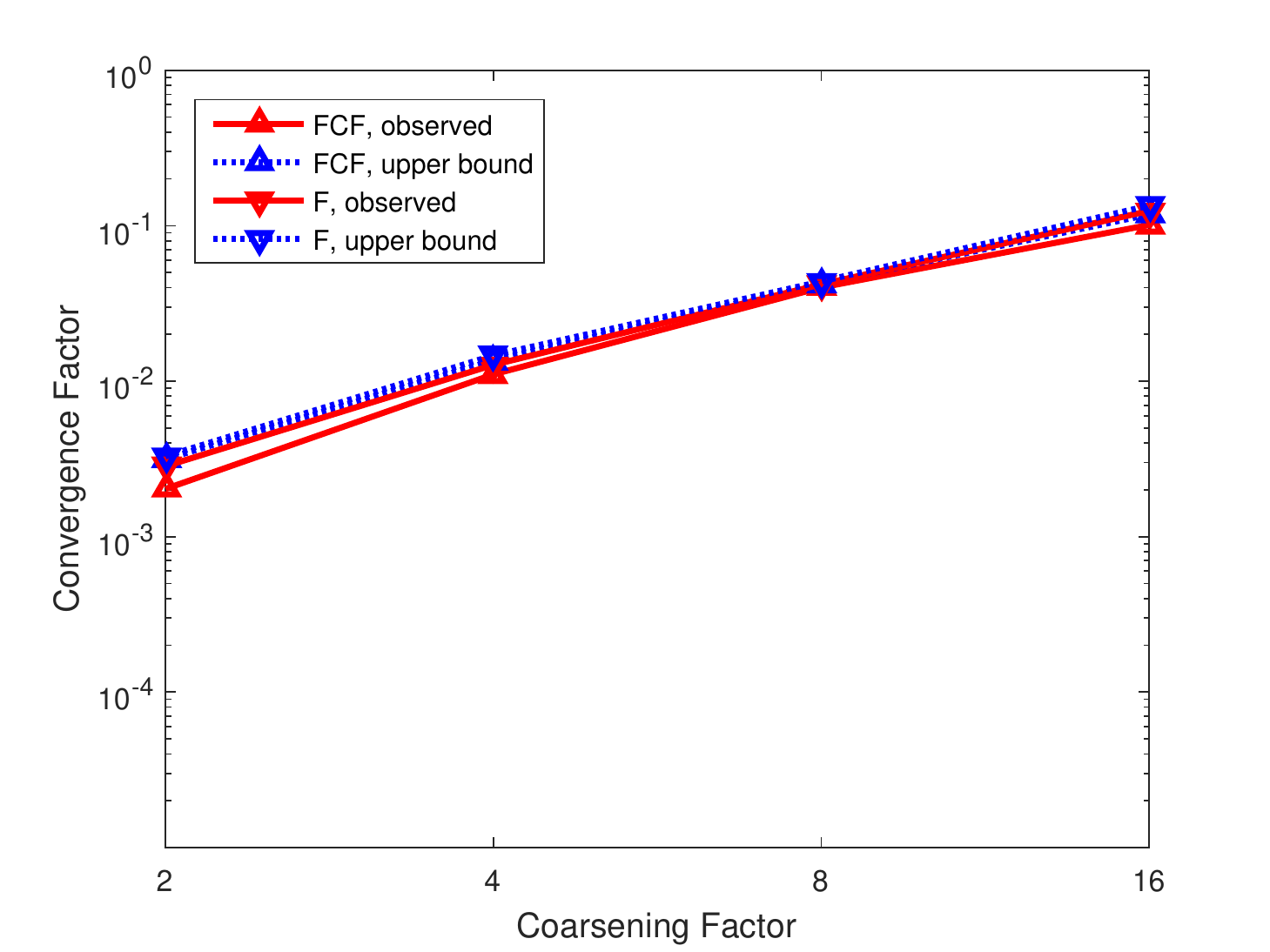}}
\hfill
\subfigure[L-shaped domain, $N=1024$]{
\label{fig-05e}
\includegraphics[scale=0.333]{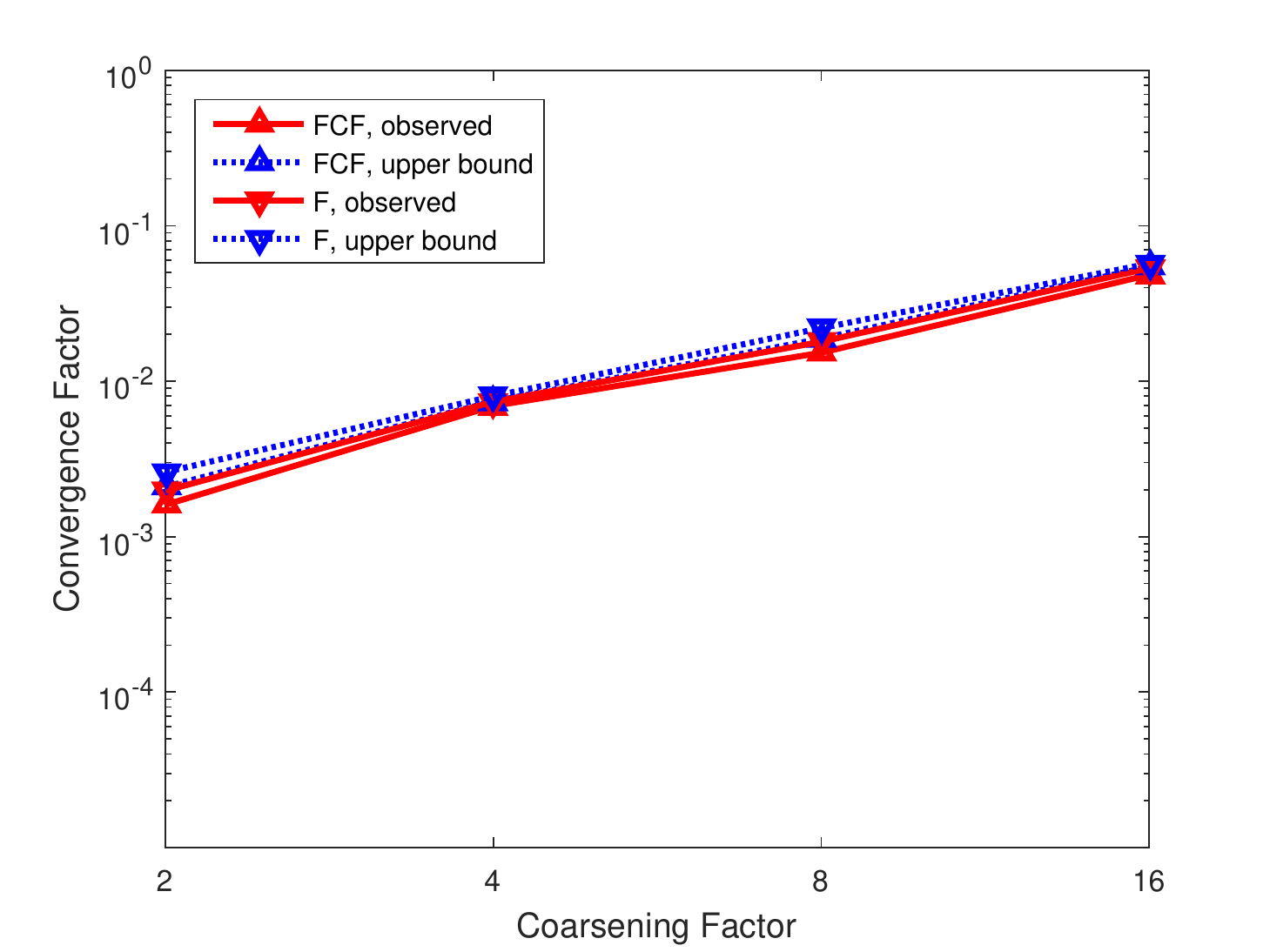}}
\hfill
\subfigure[L-shaped domain, $N=256$]{
\label{fig-05c}
\includegraphics[scale=0.333]{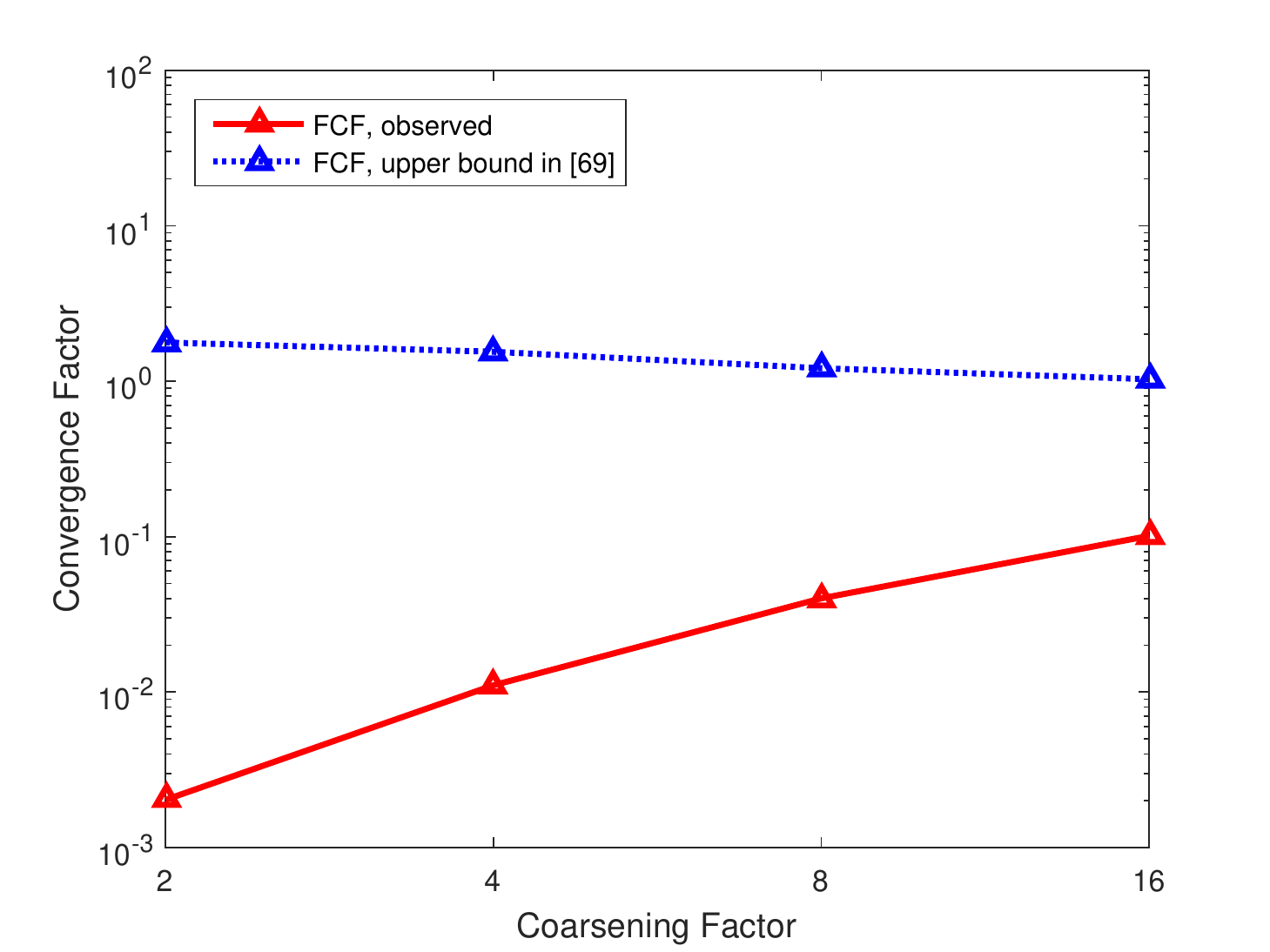}}
\caption{Graded temporal mesh: comparison on the convergence behavior of the
two-level MGRIT V(1,0)-cycle with F- and FCF-relaxation.}\label{fig-05}
\end{figure}

Fig. \ref{fig-05} reports the sharpness results in the graded temporal mesh,
where we notice again that the upper bounds on the two-level MGRIT V(1,0)-cycle
convergence described in \cite[Theorem 35]{s-03} and Theorem \ref{thm-01} are tight
(see Fig. \ref{fig-05}(a), \ref{fig-05}(b), \ref{fig-05}(d) and \ref{fig-05}(e)),
whereas the upper bound provided in \cite[Theorem 2, inequality (3.4)]{y-01} is not a good indicator,
sometimes even greater than one (see Fig. \ref{fig-05}(c) and \ref{fig-05}(f)).
Similar to the uniform case, the predicted and the largest observed convergence factors
grow with the coarsening factor but decrease with the number of time periods.
Again, FCF-relaxation yields a slightly quicker convergence.

\section{Conclusions}\label{sec5}

In this work, we provided a generalized two-level convergence theory of the MGRIT V(1,0)-cycle algorithm
with FCF-relaxation and time-dependent time-grid propagators,
and sufficient conditions for the convergence of the two-level variant in light of a new TEAP condition.
The generalization is reflected by removing the unitary diagonalization assumption, inspired by Southworth \cite{s-03}.
Numerical experiments with respect to two-dimensional unsteady fractional Laplacian problems
show consistency with the theoretical results and sharpness with our analysis.
In the future, we plan to study the multilevel convergence analysis and
extensions of this method to time fractional problems \cite{y-04,b-03,y-05,z-03,y-07,y-06}
or to support problems with adaptive time stepping \cite{h-05}
and multi-group radiation diffusion equations in inertial fusion research \cite{y-08}.

\section*{Acknowledgement}

{\color{red}The authors would like to thank the anonymous editor and reviewers for their conscientious reading and
numerous constructive comments which extremely improved the presentation of the paper.}
This work is financially supported by Science Challenge Project (TZZT2016002),
National Natural Science Foundation of China (11701197, 41874086, 11971414),
Project of Scientific Research Fund of Hunan Provincial Science and Technology Department of China (2018WK4006) and
Excellent Youth Foundation of Hunan Province of China (2018JJ1042).
We are indebted to Prof. L. Chen from University of California at Irvine for providing the software package $i$FEM \cite{c-09}.

\end{document}